\documentclass[a4paper, 10pt]{article}
\usepackage{mathtools}
\DeclareFontEncoding{FMS}{}{}
\DeclareFontSubstitution{FMS}{futm}{m}{n}
\DeclareFontEncoding{FMX}{}{}
\DeclareFontSubstitution{FMX}{futm}{m}{n}
\DeclareSymbolFont{fouriersymbols}{FMS}{futm}{m}{n}
\DeclareSymbolFont{fourierlargesymbols}{FMX}{futm}{m}{n}
\DeclareMathDelimiter{\VERT}{\mathord}{fouriersymbols}{152}{fourierlargesymbols}{147}

\usepackage[normalem]{ulem}
\usepackage{cancel}
\usepackage{bm}
\usepackage[colorlinks=true]{hyperref}
\hypersetup{urlcolor=blue,citecolor=red}
\usepackage[active]{srcltx}
\usepackage[dvips]{graphicx}
\usepackage{amssymb,amsmath,amsfonts,bbm,pifont,upgreek,bbold,accents} 
\newtheorem{theorem}{Theorem}[section]
\newtheorem{definition}{Definition}[section]
\newtheorem{proposition}{Proposition}[section]
\newtheorem{lemma}{Lemma}[section]
\newtheorem{sublemma}{Sublemma}[section]
\newtheorem{remark}{Remark}[section]
\newtheorem{notationalremark}{Notational Remark}[section]
\newtheorem{corollary}{Corollary}[section]
\newtheorem{assumption}{Assumption}[section]
\newtheorem{claim}{Claim}[section]

\newtheorem{tools}{$\negsp\negsp$}[subsection]

\newcommand\thm[1]{ \begin{theorem}\label{#1}}
\newcommand\thmtwo[2]{ \begin{theorem}[#1]\label{#2}}
\newcommand\ethm{ \end{theorem} }
\newcommand\dfn[1]{ \begin{definition}\label{#1} \rm}
\newcommand\dfntwo[2]{ \begin{definition}[#1]\label{#2} \rm}
\newcommand\edfn{ \end{definition} }
\newcommand\pro[1]{ \begin{proposition}\label{#1}}
\newcommand\protwo[2]{ \begin{proposition}[#1]\label{#2}}
\newcommand\epro{ \end{proposition} }
\newcommand\lem[1]{ \begin{lemma}\label{#1}}
\newcommand\lemtwo[2]{ \begin{lemma}[#1]\label{#2}}
\newcommand\elem{ \end{lemma} }
\newcommand\sublem[1]{ \begin{sublemma}\label{#1}}
\newcommand\sublemtwo[2]{ \begin{sublemma}[#1]\label{#2}}
\newcommand\esublem{ \end{sublemma} }
\newcommand\rem[1]{ \begin{remark}\label{#1} \rm}
\newcommand\erem{ \end{remark} }
\newcommand\notrem[1]{ \begin{notationalremark}\label{#1} \rm}
\newcommand\enotrem{ \end{notationalremark} }
\newcommand\cor[1]{ \begin{corollary}\label{#1}}
\newcommand\cortwo[2]{ \begin{corollary}[#1]\label{#2}}
\newcommand\ecor{ \end{corollary} }
\newcommand\asmp[1]{ \begin{assumption}\label{#1}}
\newcommand\asmptwo[2]{ \begin{assumption}[#1]\label{#2}}
\newcommand\easmp{ \end{assumption} }
\newcommand\clm[1]{ \begin{claim}\label{#1}}
\newcommand\eclm{ \end{claim} }
\newcommand{\proof}{\par\medskip\noindent{\bf Proof}}
\newcommand\ovl[1]{ \overline {#1} }
\newcommand\su[1]{ \frac{1}{ {#1}} }

\renewcommand{\natural}{ {\mathbf N}   }
\newcommand{{\real}}{ {\mathbf R}   }
\newcommand{{\integer}}{ {\mathbf Z}   }

\renewcommand{\d}{ {\delta}   }

\renewcommand{\l}{ {\lambda}   }

\renewcommand{\P}{ {\Pi}   }
\renewcommand{\r}{ {\rho}   }
\newcommand{\s}{ {\sigma}   }

\renewcommand{\t}{ {\tau}   }

\renewcommand{\Im}{{\, \rm Im\, }}

\newcommand{\cK}{ {\cal K} }

\newcommand{\cL}{ {\cal L} }

\newcommand{\cN}{ {\cal N} }

\newcommand{{\cJ}}{ {\cal J} }

\newcommand\ppu{{ (1) }}

\newcommand\ppj{{ (j) }}

\newcommand\ppi{{ {\rm (i)} }}

\newcommand\ppo{{ (0) }}

\usepackage{color}
\definecolor{applegreen}{rgb}{0.55, 0.71, 0.0}
\definecolor{amber(sae/ece)}{rgb}{1.0, 0.49, 0.0}
\definecolor{amethyst}{rgb}{0.6, 0.4, 0.8}

\newcommand\ii{{\rm i}}
\newcommand\cO{{\cal O}}
\newcommand\id{{\, \rm id \,}}

\newcommand\HH{{\rm H}}
\newcommand\II{{\rm I}}

\newcommand\hh{{\rm h}}

\begin{document}

%\preprint{}
\title{Non Quasi--Periodic Normal Form Theory} %Title of paper


\author{Gabriella Pinzari\footnote{Department of Mathematics, University of Padova, {\tt pinzari@math.unipd.it}}
\thanks{This paper is dedicated to Professor A. Chenciner, on the occasion of his 80th birthday, with deep admiration.}
%G.P. is indebted to G. Gallavotti for many highlighting  discussions concerning the necessity of ``estimating the distance of tori'' for the proof of Theorem A.  Theorem~\ref{stable toriREF} is part of the MSc thesis of X.L. This research is funded by the PRIN project  "New frontiers of Celestial Mechanics: theory and applications" and has been developed under the auspices of INdAM and GNFM.}}
%\email[]{gabriella.pinzari@math.unipd.it
}
%\thanks{}
%\affiliation{Universit\`a di Padova, Dipartimento di Matematica}


\date{March 15, 2023}
\maketitle %\maketitle must follow title, authors, abstract and \pacs

\begin{abstract}
We review a recent generalization of Normal Form Theory to systems (Hamiltonian ones or general ODEs) where the perturbing term is not periodic in one coordinate variable. The main difference with the standard case relies on the non uniqueness of the Normal Form and the total absence of the small divisors problem. The exposition is quite general, so as to allow extensions to the case of more non--periodic coordinates, and more  functional settings. Here, for simplicity, we work in the real--analytic class.
\end{abstract}
\tableofcontents
%\pacs{02. Mathematical methods in physics}% insert suggested PACS numbers in braces on next line

%, References should be done using the~\cite,~\ref, and \label commands
\newpage
\section{Overview}
 In this note, we review a recent generalization  of Normal Form Theory,  which we call Non Quasi--Periodic ({\sc nqp}) Normal Form. In two different forms (for Hamiltonian systems or more general  ODEs; see below), it turned out to be useful in connection with problems of celestial mechanics. The reader interested in such applications, as well as in  complete statements of the results quoted here, is referred to the original papers~\cite{pinzari2020, chenP2021}.

\vskip.1in
\noindent
Let us start with recalling that standard Normal Form Theory deals with conjugating a given Hamiltonian system, whose Hamiltonian is close--to--be--integrable:
\begin{align}\label{1}\HH(\II, \varphi)=\hh(\II)+\varepsilon f(\II, \varphi) \end{align}
where $\II=(\II_1, \ldots, \II_n)$, $\varphi=(\varphi_1, \ldots, \varphi_n)$ are couples of action--angle coordinates and $0\le \varepsilon\ll 1$ measures the size of the perturbing term, to a new Hamiltonian
\begin{align}\label{2}\HH'(\II', \varphi')=\hh'(\II')+\varepsilon'f'(\II', \varphi') \end{align}with a much smaller remainder; typically $\varepsilon'\sim e^{-\frac{1}{\varepsilon}}$.

\vskip.1in
\noindent
The origin of Normal Form Theory was undoubtedly perturbation theory (corresponding to just one step on Normal Form: see (i) below) and was seminal to {\sc kam} theory \cite{kolmogorov54, moser1962, arnold63c}. It goes back to H. Poincar\'e~\cite{Poincare:1892}; see~\cite{chencinerS1996, chenciner2015} for a deep  investigation of Poincar\'e's work. 
The first statement in the direction said above goes back by  N. N. Nekhorossev, in the late 1970s~\cite{nehorosev77}. He used it to prove that the semi--major axes  in the planetary system remain close to their initial values for times of the order $e^{\frac{1}{\varepsilon}}$. It is to be specified that in fact, in view of his application, Nekhorossev considered systems of the form
$$\HH(\II, \varphi)=\hh(\II)+\varepsilon f(\II, \varphi, p, q) $$
namely, where the perturbing term depends also on couples of ``rectangular'' coordinates $p=(p_1, \ldots, p_m)$, $q=(q_1, \ldots, q_m)$ (where by ``rectangular'' we just mean that the $q_i$'s need not to be angles) which do not appear in the unperturbed term $\hh(\II)$. A sensitive refinement of Nekhorossev's work appeared in~\cite{niederman04}.
\vskip.1in
\noindent
Typically, Normal Form Theory goes through  two steps:
\begin{itemize}
\item[(i)] constructing a  canonical transformation which allows to conjugate, symplectically,~\eqref{1} to~\eqref{2}, with $\varepsilon'=\varepsilon^2$. Typically, this  step involves {\it small divisors};
\item[(ii)] iterating (i) as many times as possible, in order to obtain an as small as possible remainder.
\end{itemize}
In the  paper by Nekhorossev,  the transformation in step (i) was obtained by means of a classical {\it generating functions}. Later, J. P\"oschel~\cite{poschel93} switched  to {\it time--one--flows} (see Definition~\ref{time one flows} below). Moreover, he clarified that step (i) is an interplay of four main points:
\begin{itemize}
\item[(i${}_{1}$)] {\it Choosing a functional setting.}
\item[(i${}_{2}$)] {\it Establishing general conditions of existence of time--one--flows  of a given Hamiltonian $\phi$, in such  functional setting}. 
\item[(i${}_{3}$)] {\it Discussing   existence and well--posedness of the solutions of the  {\it homological equation} 
  \begin{align}\label{intro: homeq2}-D\phi+f=\hh_1:\quad {\it function\ independent \ of\ \varphi}\end{align}
  with $D$ a suitable differential operator.}
\item[(i${}_{4}$)]  {\it Applying $({\rm i}_{2})$ to a solution $\phi$ determined in $({\rm i}_{3})$.}
\end{itemize}
We remark that the right hand side of 
\eqref{intro: homeq2} corresponds to the $\varepsilon$--correction to $\hh$, after one step: \begin{align}\label{h'}\HH'=\hh(\II)+\varepsilon\hh_1(\II)+\varepsilon^2 f'(\II, \varphi)\,.\end{align}

\noindent
If, one one hand, the use of time--one--flows   might sound an overly sophisticated  and somewhat renounceable tool, and  the condition
in (i${}_{2}$)  as 
an ``alter ego'' of  usual assumptions in inverse function theorems (which is the typical issue one is led to face when using generating functions) -- on the other hand, 
we shall  see that the scheme outlined in (i${}_{1}$)--(i${}_{4}$) retains such an amount of generality 
 to allow extensions  to more general dynamical systems.  Before entering in such details, let us make some remarks on the particularities of~\cite{poschel93}, which may (and will) be released later. 
 
 \vskip.1in
 \noindent
The most delicate point is (i${}_{3}$). Its main features are:
\paragraph{\it Small denominators problem\,$:$} In~\cite{poschel93} the operator $D$ in~\eqref{intro: homeq2} is given by
\begin{align}\label{oldD}D_\omega:=\omega\cdot\partial_\varphi \qquad \omega:=\partial_\II\hh\,.\end{align}
Solving~\eqref{intro: homeq2} with such $D_\omega$ requires to face the so called small denominator problem: some restriction of the domain of the  $\omega$'s (which reflects on a restriction of the domain of the $\II$'s)  is needed. Typically, one chooses the {\it Diophantine} condition. See~\cite{bounemouraF19} for a comprehensive discussion. 

 \paragraph{\it Uniqueness of Normal Form\,$:$}
Equation~\eqref{intro: homeq2} involves two unknowns: $\phi$ and $\hh_1$. In principle, one expects an infinite number of solutions. However, by very peculiar reasons of the case considered in~\cite{poschel93}, this not true. Indeed, in such case, under any of the conditions on $\omega$
mentioned in the previous item, 
  $\cK_{er}(D_\omega)$ consists of $\varphi$--independent functions, whence
 $\hh_1\in  \cK_{er}(D_\omega)$. As $D_\omega \phi\in \cK_{er}(D_\omega)^\perp$, a projection
of Equation~\eqref{intro: homeq2} on $\cK_{er}(D)$ gives $\hh_1=$ $\langle f\rangle_\varphi$
and $\phi$ is determined consequentially (see the previous item). By~\eqref{h'}, $$\hh'=\hh+\varepsilon \langle f\rangle_\varphi$$ is the only possible expression of $\hh'$, after one step. The argument obviously extends at any step, leading to uniqueness of Normal Form. In particular, there is no hope to have, in general, $\hh=\hh'$, even after one step.

 \vskip.1in
  \noindent
 In this note we consider the case when one of the  $\varphi_i$'s is {\it not} an angle. We shall see that in such a case there is much more freedom in the choice of the  Normal Form.
We discuss two examples, which we believe to be most significant among the others. One is formulated for Hamiltonian systems; the other for vector--fields. Apart for the different frameworks, the two examples rely with two different Normal Forms: In the former, after the normalization, the unperturbed part is $\varepsilon$--close to the initial one, precisely like it happens to the case in~\eqref{h'}. In the latter, we shall choose the new unperturbed part exactly the same as the initial one. In both such examples  the construction does not meet the small denominator problem. We remark that a similar problem for non--autonomous Hamiltonian systems has been considered  in~\cite{fortunatiW16}.
 
\vskip.1in
\noindent
A. We consider
the Hamiltonian
\begin{align}\label{Example1}\HH(\II, y, \varphi, x)=\hh(\II, y)+\varepsilon f(\II, y, \varphi, x) \,.\end{align}
where $(\II, \varphi)$ are action--angle coordinates; $(y, x)$ are some rectangular coordinates,  taken to be one--dimensional for simplicity. The perturbing term $f$ is not taken to be periodic in $x$. %, otherwise we fall in the case discussed in~\cite{poschel93}. 
%When $\varepsilon=0$, the motions are
%$$\left\{
%\begin{array}{lll}
%\II=\II_0\\
%y=y_0\\
%\varphi=\varphi_0+\omega_0 t\\
%x=x_0+v_0 t
%\end{array}
%\right.$$
%with $\omega_0:=\omega(\II_0, y_0):=\partial_\II\hh(\II_0, y_0)$, $v_0:=v(\II_0, y_0):=\partial_y\hh(\II_0, y_0)$.
%The only difference with respect to~\cite{poschel93} relies in (i${}_{c}$), as
% (i${}_{a}$), (i${}_{b}$) and (i${}_{d}$) can be taken precisely the same. 
 In such case, the operator $D$ in~\eqref{intro: homeq2} is
  \begin{align}\label{Lphi}D_{\omega, v}:=\omega\cdot \partial_\varphi +v\cdot \partial_x\end{align}
where
$\omega:=\partial_\II \hh$, $v:=\partial_y\hh$.
Differently from $D_\omega$ in~\eqref{oldD},  $\cK_{er}(D_{\omega, v})$ does not consists of $\varphi$--independent functions, regardless  $\omega$ is  resonant or not. 
The uniqueness argument discussed above then does not apply here.
As a matter of fact, under general conditions (see Lemma~\ref{solutionshomeq3}), Equation
$$D_{\omega, v}G=g$$
admits a solution $G$ for any choice of $g$. Taking $g=\hh_1-f$, one then has that for any $\varphi$--independent function $\hh_1$, one can find $\phi$  satisfying Equation~\eqref{intro: homeq2}. 
In~\cite{pinzari2020}, we made the choice $\hh_1=\langle f\rangle_\varphi$, in order to have a normal form which reduces to the standard  when $f$ is periodic in $x$. The Normal Form Theorem achieved, after many steps, with such choice  is discussed in Section~\ref{Theorem A} (Theorem~A).

\vskip.1in
\noindent
B. We consider the ODE: 
\begin{align}\label{oldVF}\dot x=X(x)\end{align}
where $X$ is
a $(n+1+m)$--dimensional vector--field of the form
\begin{align}\label{perturbed}X({\rm I}, y, \varphi)=N({\rm I}, y)+\varepsilon P({\rm I}, y, \varphi)\quad x=({\rm I}, y, \varphi)\in V\times{\mathbb T}^n\quad V\subset {\mathbb R}^m\times{\mathbb R} \end{align}
where
$$N({\rm I}, y) =
\left(
\begin{array}{c}
0\\
v({\rm I}, y)\\
\omega({\rm I}, y)
\end{array}
\right)%=
%v({\rm I}, y)\partial_y+\omega({\rm I}, y) \partial_\varphi 
\, \qquad P({\rm I}, y, \varphi )=\left(
\begin{array}{c}
P_1({\rm I}, y, \varphi )\\
P_2({\rm I}, y, \varphi )\\
P_3({\rm I}, y, \varphi )
\end{array}
\right)\,. $$
Similarly as in the previous cases, the general problem is to find a diffeomorphism  $$\Phi:\quad z\in V\times{\mathbb T}^n \to x=\Phi(z)\in U\times{\mathbb T}^n\qquad U\,, \ V\subset {\mathbb R}^m\times {\mathbb R}$$ which transforms equation~\eqref{oldVF}
to
\begin{align}\label{newVF}\dot z=Z(z)
\end{align}
with
$$Z({\rm I}, y, \varphi )=N'({\rm I}, y)+\varepsilon' P'({\rm I}, y, \varphi) $$
and $0<\varepsilon'\ll \varepsilon$.
Such a transformation will be searched in the form of a the time--one--flow
of a vector--field $Y$. The {\it homological equation} which allows to make one step, namely, with $\varepsilon'=\varepsilon^2$, is now given by (see Section~\ref{Theorem B})
$$-D_NY +P=N_1\qquad \textrm{vector--field independent of $\varphi$}$$
where
$D_NY=[N, Y]$
are the Lie brackets of $N$ and $Y$ (see Definition~\ref{Vector--field time--one  flow}).
Again,
$N_1$ corresponds to the $\varepsilon$--correction of $N$:
$$N'=N+\varepsilon N_1\,.$$
 We shall show that we have a phenomenon similar to the one discussed in example A: 
under general conditions (see Proposition~\ref{homeq1}), Equation
$$D_NY=Z$$
admits a solution for any choice of $Z$, whence $N_1$ can be arbitrarily chosen.  In~\cite{chenP2021} we took $N_1=0$, which allows to have $N=N'$. The Normal Form Theorem achieved, after many steps, with such choice  is discussed in Section~\ref{Theorem B} (Theorem~B).

\vskip.1in
\noindent
We conclude this overview with saying that the purpose of this note is  to give a quite general exposition, so as to allow generalizations, if needed. As an example, a version of Theorem~B for smooth vector--fields  has been stated and used in~\cite{chenP2021}, in connection with a  procedure of regularization of collisions in the three--body problem, after which  the  holomorphy and Hamiltonian nature of the system are destroyed.

\section{Theorem~A}\label{Theorem A}

  We allow $\HH$ slightly more general than~\eqref{Example1}, namely,
  \begin{align}\label{h0f0}\HH(\II, \varphi,  p ,  q , y, x)=\hh( I , J ( p ,  q ), y)+f( \II ,  \varphi ,  p ,  q , y, x)\end{align}
%$$\hh(y, I , J )= J ( I , y)+{\mathbf P}(y, I , J )\ ,\qquad {\mathbf P}(y, I ,0)\equiv 0$$
where $p=(p_1,\cdots, p_m)$, $q=(q_1,\cdots, q_m)$ and
$J ( p ,  q )=(p_1q_1,\cdots, p_mq_m)$. 

\subsection*{Functional Setting and Statement}

We assume that $\HH$ in~\eqref{h0f0} is holomorphic\footnote{As usual, if $A\subset {\mathbb R}$ and $r${,$s>0$},  the symbols $A_r${, ${\mathbb T}_s$} denote the complex $r$, {$s$}--neighbourhoods of $A${,${\mathbb T}$:}
\begin{align}A_r:=\bigcup_{x\in A}B_r(x)\,,\qquad {{\mathbb T}_s:=\big\{\varphi =\varphi _1+{\rm i}\varphi _2:\ \varphi _1\in {\mathbb T}\,,\ \varphi _2\in {\mathbb R}\,,\ |\varphi _2|<s\big\}\,,}\end{align} with $B_r(x)$ being the complex ball centred at $x$ with radius $r$.} in
$${\mathbb P}_{\r, s, \delta, r, \xi}={\mathbb I}_\r  \times {{\mathbb T}}^n_s\times {\mathbb B}_\delta\times {\mathbb Y}_r\times  {\mathbb X}_\xi\,.$$

\noindent
We denote as ${\cal O}_{\r, s, \d, r, \xi}$ the set of complex holomorphic functions $$\phi:\quad {\mathbb P}_{\r, s,  \d,  r, \xi}\to {{\mathbb C}}$$ We equip ${\cal O}_{\r, s, \d, r, \xi}$ with 
 the norm
$$\|\phi\|_{\r, s, \d, r, \xi}:=\sum_{k,h,j}\|\phi_{khj}\|_{\r, r,\xi}e^{s|k|}\d^{h+j}$$
where $\phi_{khj}( I , y,  x)$ are the coefficients of the Taylor--Fourier expansion\footnote{We denote as ${\mathbf x}^h:=x_1^{h_1}\cdots x_n^{h_n}$, where ${\mathbf x}=(x_1, \cdots, x_n)\in {\real}^n$ and $h=(h_1, \cdots, h_n)\in {\natural}^n$.}
\begin{align}\label{expphi}\phi=\sum_{k,h,j}\phi_{khj}( I ,y, x)e^{\ii k s} p ^h  q ^j\ ,\quad \|\phi\|_{\r, r,\xi}:=\sup_{{\mathbb I}_\r\times{\mathbb Y}_r\times  {\mathbb X}_\xi}|\phi( I , y, x)|\ .\end{align}
If $\phi$ is independent of $x$, we simply write
$\|\phi\|_{\r, r}$ for $\|\phi\|_{\r, r,\xi}$.

%\noindent
%For a given vector--valued  function $\underline\phi=(\phi_1,\cdots, \phi_k)\in {\cal O}_{\r, s, \d, r, \xi}^k$, we let
%$$\|\underline\phi\|_{\r, s, \d, r, \xi}:=\sum_{i=1}^k \|\phi_i\|_{\r, s, \d, r, \xi}\ .$$
\noindent
We decompose
$${\cal O}_{\r, s, \d, r, \xi}={\cal Z}_{\r, s, \d, r, \xi}\oplus {\cal N}_{\r, s, \d, r, \xi}\,.$$via the identities (with $\phi$ as in~\eqref{expphi})

\begin{eqnarray*}&&\widetilde\phi:=\P_{{\cal Z}_{\r, s, \d, r, \xi}}\phi=\sum_{k,h,j:\atop (k,h-j)\ne (0,0)}\phi_{khj}( I , y,  x)e^{\ii k s} p ^h  q ^j\nonumber\\
&&\ovl\phi:=\P_{{\cal N}_{\r, s, \d, r, \xi}}\phi=\sum_{h:}\phi_{0hh}( I , y,  x) (pq)^h 
\,.\end{eqnarray*}

\noindent
Finally, given $\HH$ as in~\eqref{h0f0}, we let
$$\omega:=\partial_\II\hh\,,\qquad \omega':=\partial_J\hh\,,\qquad v:=\partial_y\hh\,.$$

\noindent
We shall prove the following

\vskip.1in
\noindent
{\bf Theorem~A}\ {\it Let the Hamiltonian $\HH$  in~\eqref{h0f0} be holomorphic in ${\mathbb P}_{\r, s, \delta, r, \xi}$.
Put ${\rm d}:=\min\big\{\r s, r\xi, {\d}^2\big\}$, ${\cal X}:=\sup\big\{|x|:\ x\in {\mathbb X}_\xi\big\}$.
There exists a number ${\rm c}\ge 1$, depending only on $n$, $m$, such that, for any $p\in {\mathbb N}$ such that the following inequalities are satisfied
\begin{align}\label{normal form assumptions}
4p{\cal X}\left\|\Im\frac{\omega }{v }\right\|_{\r, r}<
s
\ ,\quad 
4p{\cal X}\left\|\frac{{\omega'} }{v }\right\|_{\r, r}<
1\,,\quad
  {{\rm c}p\frac{{\cal X}}{{\rm d}}   % \|\frac{1}{v }\|_{\r, r}
 \left\|{f}\right\|_{\r, s,  r, \xi}\left\|\frac{1}{v }\right\|_{\r, s,  r, \xi} <1}%\nonumber\\
  %&& {\rm c}^{\red2}_{n,m}p^{\red2}\frac{{\cal X}^{2}}{{\rm d}^{2}}  { \left\|\frac{1}{v }\right\|^2_{\r, r} 
%\left\|% \frac
%{f}%{v }
%\right\|_{\r, s,  r, \xi}}\left\|% \frac
%\widetilde{f}%{v }
%\right\|_{\r, s,  r, \xi}<1 
\end{align}
one can find a real--analytic transformation $\Psi_*$
which carries $\HH$ to
$$\HH_*=\hh+g_*+f_*$$
where
$g_*$ is $\varphi$--independent and, moreover, the following inequalities hold 
\begin{align}\label{thesis}
&&\|g_*-\ovl f\|_{1/3 (\r, s,  r, \xi)}\le 162\,{\rm c} \frac{{\cal X}}{\rm d}\left\|%\frac{1}{v }\|_{\r, r} \|
{\frac{\widetilde f}{v }}\right\|_{\r, s,  r, \xi}\| f\|_{\r, s,  r, \xi}\nonumber\\
&&    \|f_*\|_{1/3 (\r, s, r, \xi)}\le {2^{-(p+1)}} \|f\|_{\r, s, r, \xi}\ .\end{align}
%The transformation $\Psi_*$ can be obtained as a composition of time--one Hamiltonian flows, and satisfies the following.  If
%\begin{align}( \II ,  \varphi ,  p ,  q , y, x):=\Psi_*( I _*,  \varphi _*,  p _*,  q _*, {\rm R}_*, {\rm r}_*)\end{align} the following uniform bounds hold:
%\begin{align}\label{phi close to id***}
%&&{\rm d}\max\Big\{\frac{| I - I _*|}{\r},\ \frac{| \varphi  - \varphi  _*|}{s},\ 
%\frac{| p - p _*|}{\d},\ \frac{ | q - q _*|}{\d},\ \frac{|y-y_*|}{r},\frac{ |x-x_*|}{\xi} \Big\}\nonumber\\
%&&\le\max\Big\{s| I - I _*|,\ \r| \varphi  - \varphi  _*|,\ 
%\d | p - p _*|,\ \d | q - q _*|,\ \xi |y-y_*|,\ r |x-x_*| \Big\}\nonumber\\
%&&\leq
% {{19\, {\cal X}   \left \| \frac{ f}{v }\right\|_{\r, s,  r, \xi}
%}} \ .
%\end{align}
}

\subsection*{The Lie Series}
\begin{definition}\label{time one flows}\rm If
\begin{align}\label{J}J=\left(
\begin{array}{lc}
\mathbb{0}&-\mathbb{1}_n\\
\mathbb{1}_n&\mathbb{0}
\end{array}
\right)\end{align}
is the symplectic unit and $f=f(y)$, $g=g(y)$ are any two functions $ C^{\infty}(U)$, with $U\subset{\mathbb R}^{2n}$, we denote as
$$\{f, g\}:=J\partial_y f\cdot \partial_y g=\partial_p f\cdot \partial_q g-\partial_q f\cdot \partial_p g \qquad y=(p, q)$$ their {\it Poisson parentheses}.
We denote as
$${\cal L}_\phi :=\big\{\phi\,,\cdot\big\}:\qquad C^{\infty}(U)\to C^{\infty}(U)$$
the  {\it Hamiltonian Lie operator} (or, simply,  {\it Lie operator}, when there is no risk of confusion) on $C^{\infty}(U)$, with
a fixed $\phi\in C^\infty (U)$. The map 
\begin{align}\label{HamLieSeries}e^{\cL_\phi}:=\sum_{k=0}^\infty\frac{{\cal L}_\phi^k}{k!}\end{align}will be called the  {\it Hamiltonian Lie series} (respectively, {\it Lie series}), generated by $\phi$. 
\end{definition}
The importance of the Lie series  in the following

\begin{proposition}\label{prop: time one map}
Assume that
 the Lie series $e^{{\cal L}_\phi} $ converges uniformly. Then  the time--one map of $\phi$, $\Phi^\phi_1$,
 carries the Hamiltonian $\HH$ to  
$\HH_1=e^{{\cal L}_\phi} \HH$.
\end{proposition}

\noindent
The proof of Proposition~\ref{prop: time one map} can be obtained as an application of a more general result  (Proposition~\ref{prop: Lie operator} below) to the Hamiltonian vector--fields $X=J\partial \HH$, $Y=J\partial \phi$. However, for sake of self--consistency, we propose an independent proof in Appendix~\ref{app: time one flow}.

\subsection*{Convergence of the Lie Series}
The convergence of the series~\eqref{HamLieSeries} in the real--analytic context is well known since~\cite{poschel93}:

\begin{lemma}[\cite{poschel93}]\label{base lemma}
There exists an integer number $\ovl{\rm c}$ such that, for any $\phi\in {\cal O}_{\r, s, \d, r, \xi}$ and any $r'<r$, $s'<s$, $\r'<\r$, $\xi'<\xi$, $\d'<\d$ such that
$$\frac{\ovl{\rm c}\|\phi\|_{\r, s, \d, r, \xi}}{d}<1\qquad d:=\min\big\{\r'\s', r'\xi', {\d'}^2\big\} $$
then  the series~\eqref{HamLieSeries}
converges  uniformly to an operator
$$ e^{{\cal L}_\phi}:\quad  {\cal O}_{\r, s, \d, r, \xi}\to \cO_{\r-\r', s-s', \d-\d', r-r', \xi-\xi'}$$
and its tails
\begin{align}\label{queue}\Phi_h:=\sum_{k=h}^{\infty}\frac{{\cal L}_\phi^k}{k!}\ .\end{align}
verify\footnote{In \cite{poschel93} it is proved that \begin{align}\label{geometric series}\|\cL^k_\phi[g]\|_{\r-\r', s-s', \d-\d', r-r', \xi-\xi'}\le k!\big(\frac{\ovl{\rm c}\|\phi\|_{\r, s, \d, r, \xi}}{d}\big)^k\|g\|_{\r, s, \d, r, \xi}\end{align}
which immediately implies~\eqref{h power}.}
\begin{align}\label{h power}\|\Phi_hg\|_{\r-\r', s-s', \d-\d', r-r', \xi-\xi'}\le \frac{\big(\frac{\ovl{\rm c}\|\phi\|_{\r, s, \d, r, \xi}}{d}\big)^h}{1-\frac{\ovl{\rm c}\|\phi\|_{\r, s, \d, r, \xi}}{d}}\|g\|_{\r, s, \d, r, \xi}\qquad \forall g\in {\cal O}_{\r, s, \d, r, \xi}\ .\end{align}
\end{lemma}

\subsection*{The Homological Equation}
Define
\begin{align}\label{Lphi}{\cal L}_\phi \hh:=\{\phi, \hh\}=-D_{\omega, \omega', v}\phi\end{align}
where, if $\phi$ is as~\eqref{expphi}, 
then
$$
D_{\omega, \omega', v}\phi:=\sum_{(k,h,j)} \big(v \partial_x +(h-j)\cdot{\omega'} +\ii k\cdot\omega \big)\phi_{khj}( I ,y, x)e^{\ii k\cdot \varphi } p ^h  q ^j\,.$$

\begin{lemma}\label{solutionshomeq3}
The function
\begin{align}\label{G}G(\II, y, \varphi, x):= \frac{1}{v }\int_0^x g\left(\II, \varphi+\frac{\omega }{v }(\t-x), p e^{\frac{{\omega'} }{v }(\t-x)}, q e^{-\frac{{\omega'} }{v }(\t-x)}, y, \t\right)d\t\end{align}
solves
\begin{align}\label{homeq3}D_{\omega, \omega', v}G=g\end{align}
Any other solution of~\eqref{homeq3} differs from $G$ for a function of the form 
\begin{align}\label{F}F=F\left(\II, \varphi-\frac{x }{v } \omega, p e^{\frac{-{\omega'} }{v }x}, q e^{\frac{{\omega'} }{v }x}, y\right)\,.\end{align}
\end{lemma}
\proof\, That~\eqref{G} solves~\eqref{homeq3} is an easy check. On the other hand, the functions~\eqref{F} are precisely the ones in the kernel of $D_{\omega, \omega', v}$. $\quad \square$

\subsection*{The Step Lemma}
Let
  \begin{align}\label{h0f0g0}\HH(\II, \varphi,  p ,  q , y, x)=\hh( I , J ( p ,  q ), y)+g( I , J ( p ,  q ), y)+f( \II ,  \varphi ,  p ,  q , y, x)\end{align}
\begin{lemma}\label{iterative lemma}
There exists a number $\widetilde{\rm c}>1$ such that the following holds. For any choice  of positive numbers  $r'$, $\r'$, $s'$, $\xi'$. $\d'$ satisfying
%\begin{align}
%&&r'<r\ ,\quad \r'<\r\ ,\quad \xi'<\xi\\
%&& s'<s\ ,\quad \d'<\d\ ,\quad {\cal X}\left\|\frac{\omega }{v }\right\|_{\r, r}<s-s'\ ,\quad
%{\cal X}\left\|\frac{{\omega'} }{v }\right\|_{\r, r}<
%\log\frac{\d}{\d'} 
%\end{align}
\begin{eqnarray}\label{ineq1}
&&{2\r'<\r\ ,\quad 2r'<r\ ,\quad  2\xi'<\xi}\\
\label{ineq2}
&&2s'<s\ ,\quad 2\d'<\d\ ,\quad { {\cal X}\left\|\Im\frac{\omega }{v }\right\|_{\r, r}<s-2s'\ ,\quad
{\cal X}\left\|\frac{{\omega'} }{v }\right\|_{\r, r}<
\log\frac{\d}{2\d'} }
\end{eqnarray}
and and provided that the  following inequality holds true
\begin{align}\label{smallness}
 \widetilde{\rm c}\frac{{\cal X}}{d}   % \|\frac{1}{v }\|_{\r, r}
 \left \|\frac{\widetilde f}{v }\right\|_{\r, s, \d, r, \xi} <1\qquad d:=\min\big\{\r'\s', r'\xi', {\d'}^2\big\}
\end{align}
 one can find an operator $$\Phi:\quad {\cal O}_{\r, s, \d, r, \xi}\to \cO_{\r_+, s_+, \d_+, r_+, \xi_+} $$ with
%  $$r_+:=r-r'\ ,\quad \r_+:=\r-\r'\ ,\quad \xi_+:=\xi-\xi'\ ,\quad s_+:=s-s'-{\cal X}\left\|\frac{\omega }{v }\right\|_{\r, r}\ ,\quad \d_+:=\d e^{-{\cal X}\left\|\frac{{\omega'} }{v }\right\|_{\r, r}}-\d'$$
\begin{eqnarray*}
&&r_+:=r-2r'\ ,\quad \r_+:=\r-2\r'\ ,\quad \xi_+:=\xi-2\xi'\ ,\quad s_+:=s-2s'-{\cal X}\left\|\Im\frac{\omega }{v }\right\|_{\r, r}\nonumber\\
&& \d_+:=\d e^{-{\cal X}\left\|\frac{{\omega'} }{v }\right\|_{\r, r}}-2\d'
\end{eqnarray*}
which carries the Hamiltonian  $\HH$ in~\eqref{h0f0g0} to
$$\HH_+:=\Phi[\HH]=\hh+g+\ovl f+f_+ $$
where 
\begin{align}\label{bound}\|f_+\|_{r_+,\r_+, \xi_+, s_+,\d_+}\le \widetilde{\rm c}\left(\frac{{\cal X}}{d}\left \|\frac{\widetilde f}{v }\right\|_{\r, s, \d, r, \xi}\| f\|_{\r, s, \d, r, \xi}+ \|\{\phi, g\}\|_{\r_1-\r', s_1-s', \d_1-\d', r_1-r', \xi_1-\xi'}\right)
\end{align}
with $$\r_1:=\r\ ,\quad s_1:=s-{\cal X}\left\|\Im\frac{\omega }{v }\right\|_{\r, r}\ ,\quad \d_1:=\d e^{-{\cal X}\left\|\frac{{\omega'} }{v }\right\|_{\r, r}}\ ,\quad r_1:=r\ ,\quad  \xi_1:=\xi $$for a suitable $\phi\in \cO_{\r_1, s_1, \d_1, r_1, \xi_1}$ verifying \begin{align}\label{bound on phi}\|\phi\|_{\r_1, s_1, \d_1, r_1, \xi_1 }\le {\cal X}   %\|\frac{1}{v }\|_{\r, r} \|\widetilde f\|_{\r, s, \d, r, \xi}
\left \|\frac{\widetilde f}{v }\right\|_{\r, s, \d, r, \xi}
\ .\end{align}\end{lemma}

\proof\, 
Define
\begin{eqnarray*}
&&\ovl r=r%=:r_1
\ ,\quad\ovl\r=\r%=: \r_1
\ ,\quad \ovl\xi=\xi%=:\xi_1
\ ,\quad \ovl s=s-{\cal X}\left\|\Im\frac{\omega }{v }\right\|_{\r, r}%=:s_1
%\nonumber\\
%&&
\ ,\quad
 \ovl\d=\d e^{-{\cal X}\left\|\frac{{\omega'} }{v }\right\|_{\r, r}}%=:\d_1
\end{eqnarray*}
In view of the two last but one inequality in~\eqref{ineq2}, we have
$0<\ovl s< s$.
%and assume that
%\begin{align}\label{ovl param1}{\cal X}\left\|\Im\frac{\omega }{v }\right\|_{ \r, r} \le s-\ovl s\ ,\qquad {\cal X}\left\|\frac{{\omega'} }{v }\right\|_{ \r, r} \le \log\frac{\d}{\ovl \d}\ .\end{align}
The formal application of $\Phi=e^{\cL_\phi}$ to the Hamiltonian~\eqref{h0f0g0} yields:
\begin{align}\label{f1}
e^{\cL_\phi} \HH&=&e^{\cL_\phi} (\hh+g+
f)=\hh+g+
{\cal L}_\phi \hh+f+\Phi_2(\hh)+\Phi_1(g)
+\Phi_1(f)
\end{align}
where the $\Phi_h$'s are the tails of $e^{\cL_\phi}$, defined in~\eqref{queue}.
We choose $\phi$ so that
\begin{align}\label{homeq4}{\cal L}_\phi \hh+f=\ovl f\,.\end{align}
By~\eqref{Lphi} and Lemma~\ref{solutionshomeq3}, the function
$$\phi(\II, y, \varphi, x):= \frac{1}{v }\int_0^x \widetilde f\left(\II, \varphi+\frac{\omega }{v }(\t-x), p e^{\frac{{\omega'} }{v }(\t-x)}, q e^{-\frac{{\omega'} }{v }(\t-x)}, y, \t\right)d\t$$
solves~\eqref{homeq4}.
Then we have
$$\|\phi_{khj}\|_{\ovl\r,\ovl r,\ovl\xi}\le  %\|\frac{1}{v }\|_{\ovl\r, \ovl r}
 \left\|\frac{f_{khj}}{v }\right\|_{\ovl\r,\ovl r,\ovl\xi}\left\|\int_0^x |e^{-\frac{\l_{khj}}{v }\t}|\right\|_{\ovl\r,\ovl r,\ovl\xi}d\t\le{\cal X}   %\|\frac{1}{v }\|_{\ovl\r, \ovl r} 
  \left\|\frac{f_{khj}}{v }\right\|_{\ovl\r,\ovl r,\ovl\xi} e^{{\cal X}\left\|\frac{\l_{khj}}{v }\right\|_{\ovl\r, \ovl r}}\, . $$
%We take, in particular, $\ovl r$, $\ovl\r$, ${\cal X}$ as above and, moreover,
%$\ovl s$, $\ovl\d$ such that
% $$\ovl s+{\cal X}\left\|\frac{\omega }{v }\right\|_{\ovl\r, \ovl r} \le s\ ,\qquad \ovl \d e^{{\cal X}\left\|\frac{{\omega'} }{v }\right\|_{\ovl\r, \ovl r} }\le \d\ .$$
Since
$$\|\frac{\l_{khj}}{v }\|_{\ovl\r, \ovl r} \le (h+j)\left\|\frac{{\omega'} }{v }\right\|_{\ovl\r, \ovl r} +|k|\left\|\Im\frac{\omega }{v }\right\|_{\ovl\r, \ovl r} 
$$
we have
$$\|\phi_{khj}\|_{\ovl\r,\ovl r,\ovl\xi}\le {\cal X}   %\|\frac{1}{v }\|_{\ovl\r, \ovl r}  
\left\|\frac{\widetilde f_{khj}}{v }\right\|_{\ovl\r,\ovl r,\ovl\xi}
e^{(h+j){\cal X}\left\|\frac{{\omega'} }{v }\right\|_{\ovl\r, \ovl r} +|k|{\cal X}\left\|\Im\frac{\omega }{v }\right\|_{\ovl\r, \ovl r} }\, .
$$
which yields (after multiplying by $e^{|k|\ovl s}(\ovl\d)^{j+h}$ and summing over $k$, $j$, $h$ with $(k,h-k)\ne (0,0)$) to
$$\|\phi\|_{\ovl\r,\ovl r,\ovl\xi, \ovl s,\ovl\d}\le{\cal X}  %\|\frac{1}{v }\|_{\ovl r,\ovl \r} 
\left\|\frac{\widetilde f}{v }\right\|_{\ovl r,\ovl \r,\ovl \xi,\ovl s+{\cal X}\left\|\Im\frac{\omega }{v }\right\|_{\ovl\r, \ovl r} ,\ovl \d e^{{\cal X}\left\|\frac{{\omega'} }{v }\right\|_{\ovl\r, \ovl r} }}=
{\cal X}  %\|\frac{1}{v }\|_{\ovl r,\ovl \r} 
\left\|\frac{\widetilde f}{v }\right\|_{r, \r, \xi, s ,\delta}
\ . $$
An application of Lemma~\ref{base lemma},with $r$, $\r$, $\xi$, $s$, $\d$ replaced by $\ovl r-r'$, $\ovl r-\r'$, $\ovl \xi-\xi'$, $\ovl s-s'$, $\ovl \d-\d'$,  concludes  with a suitable choice of $\widetilde{\rm c}>\ovl{\rm c}$ and (by~\eqref{f1}) $$ f_+:=\Phi_2(\hh)+\Phi_1(g)+\Phi_1(f)\ . $$
Observe that the bound~\eqref{bound} follows from~\eqref{h power},~\eqref{geometric series} and the identities
$$\Phi_2[\hh]=\sum_{j=2}^\infty \frac{\cL^j_\phi(\hh)}{j!}=\sum_{j=1}^\infty \frac{\cL^{j+1}_\phi(\hh)}{(j+1)!}=-\sum_{j=1}^\infty \frac{\cL^{j}_\phi(\widetilde f)}{(j+1)!}
$$
$$\Phi_1[g]=\sum_{j=1}^\infty \frac{\cL^j_\phi(g)}{j!}=\sum_{j=0}^\infty \frac{\cL^{j+1}_\phi(g)}{(j+1)!}=\sum_{j=0}^\infty \frac{\cL^{j}_\phi(g_1)}{(j+1)!} $$
with $g_1:=\cL_\phi(g)=\{\phi, g\}$. $\qquad\square$

\subsection*{Iterations (Proof of Theorem~A)}
The proof of the  Theorem~A goes through iterate applications of Lemma~\ref{iterative lemma}. We premise the following

\begin{remark}\label{stronger iterative lemma}\rm
Replacing conditions in~\eqref{ineq2} with the stronger ones\footnote{The three first inequalities in~\eqref{new cond} are immediately seen to be stronger that the corresponding three first inequalities  in~\eqref{ineq2}.
On the other hand, rewriting the second inequality in~\eqref{new cond} as
$\frac{\d'}{\d}<1-\frac{2\d'}{\d}$
and using the inequality (which holds for all $x\ge 1$) $\log x\ge 1-\frac{1}{x}$  with $x=\frac{\d}{2\d'}$, we have  also $\frac{\d'}{\d}<\log\frac{\d}{2\d'} $. 
}
%\beqno2s'<s\ ,\quad 2\d'<\d\ ,\quad {\cal X}\left\|\frac{\omega }{v }\right\|_{\r, r}<s'\ ,\quad
%{\cal X}\left\|\frac{{\omega'} }{v }\right\|_{\r, r}<
%\frac{\d'}{\d} \eeqno

\begin{align}\label{new cond}
{3s'<s\ ,\quad 3\d'<\d\ ,\quad {\cal X}\left\|\Im\frac{\omega }{v }\right\|_{\r, r}<s'\ ,\quad
{\cal X}\left\|\frac{{\omega'} }{v }\right\|_{\r, r}<
\frac{\d'}{\d}} \end{align}

\noindent
(and keeping~\eqref{ineq1},~\eqref{smallness} unvaried)  one can take, for $s_+$, $\d_+$,  $s_1$, $\d_1$ the simpler expressions
%$$s_{+\rm new}=s-2s'\ ,\quad \d_{+\rm new}=\d-2\d'\ ,\quad s_{1\rm new}:=s-s'\ ,\quad \d_{1\rm new}=\d-\d'\ $$

$${s_{+\rm new}=s-3s'\ ,\quad \d_{+\rm new}=\d-3\d'\ ,\quad s_{1\rm new}:=s-s'\ ,\quad \d_{1\rm new}=\d-\d'\ }$$
(while keeping $r_+$, $\r_+$, $\xi_+$, $r_1$, $\r_1$, $\xi_1$ unvaried).
Indeed, since $1-e^{-x}\le x$ for all $x$,
$$\d_1=\d e^{-{\cal X}\left\|\frac{{\omega'} }{v }\right\|_{\r, r}} = \d-\d(1- e^{-{\cal X}\left\|\frac{{\omega'} }{v }\right\|_{\r, r}})\ge \d-\d{\cal X}\left\|\frac{{\omega'} }{v }\right\|_{\r, r}\ge \d-\d' =\d_{1\rm new}\ . $$
This also implies $\xi_+=\d_1-\d'\ge \d-2\d'=\xi_{+\rm new}$. That $s_+\ge s_{+\rm new}$, $s_1\ge s_{1\rm new}$ is even more immediate.\end{remark}

\noindent
We apply Lemma~\ref{iterative lemma} with
$${2}\r'=\frac{\r}{3}\ ,\quad {3}s'=\frac{s}{3}\ ,\quad {3}\d'=\frac{\d}{3}\ ,\quad{2}r'=\frac{r}{3}\ ,\quad {2}\xi'=\frac{\xi}{3}\ ,\quad g\equiv 0\, . $$
We make use of the stronger formulation described in Remark~\ref{stronger iterative lemma}.  Conditions in~\eqref{ineq1} and the two former conditions in~\eqref{new cond} are trivially true. The two latter inequalities in 
\eqref{new cond} reduce to 
$${\cal X}\left\|\Im\frac{\omega }{v }\right\|_{\r, r}<\frac{s}{9}\ ,\quad
{\cal X}\left\|\frac{{\omega'} }{v }\right\|_{\r, r}<
\frac{1}{9} $$and they are certainly satisfied by assumption~\eqref{normal form assumptions}, for {$p>2$}. Since  
$$ d=\min\{ \r' s', r'\xi', {\d'}^2 \}=\min\{ \r s/{36}, r\xi/{54}, {\d}^2/{81} \}\ge\frac{1}{{81}}
\min\{ \r s, r\xi, {\d}^2 \}=\frac{\rm d}{{81}}
$$
we have that condition~\eqref{smallness} is certainly implied by the last inequality in~\eqref{normal form assumptions}, once one chooses ${\rm c}>{81} \widetilde{\rm c}$.
By Lemma~\ref{iterative lemma}, it is then possible to conjugate $\HH$ to
$$\HH_1=\hh_0+\ovl f+f_1$$with $f_1\in \cO_{\r^\ppu,s^\ppu,\d^\ppu, r^\ppu,\xi^\ppu,}$, where $(\r^\ppu,s^\ppu,\d^\ppu, r^\ppu,\xi^\ppu,):=2/3 (\r, s, \d, r, \xi)$ and
\begin{align}\label{f1***}\|f_1\|_{\r^\ppu,s^\ppu,\d^\ppu, r^\ppu,\xi^\ppu,}\le {81}\widetilde{\rm c} \frac{{\cal X}}{\rm d}  % \|\frac{1}{v }\|_{\r, r} 
\left\|\frac{\widetilde f}{v }\right\|_{\r, s, \d, r, \xi}\| f\|_{\r, s, \d, r, \xi}\le \frac{\| f\|_{\r, s, \d, r, \xi}}{2}\end{align}
since ${\rm c}\ge {162} \widetilde{\rm c}$ and $p\ge1$. Now we aim to apply Lemma~\ref{iterative lemma} $p$ times, again as described  in Remark~\ref{stronger iterative lemma}, each time with parameters
$$\r_j'=\frac{\r}{6p}\ ,\quad s_j'=\frac{s}{{9}p}\ ,\quad \d_j'=\frac{\d}{{9}p}\ ,\quad r_j'=\frac{r}{6p}\ ,\quad \xi_j'=\frac{\xi}{6p}\, .$$
Therefore we find, at each step
\begin{eqnarray*}
&&\r^{(j+1)}:=\r^\ppu-j\frac{\r}{3p}\ ,\quad s^{(j+1)}:=s^\ppu-j\frac{s}{3p}\ ,\quad \d^{(j+1)}:=\d^\ppu-j\frac{\d}{3p}\nonumber\\
&&r^{(j+1)}:=r^\ppu-j\frac{r}{3p}\ ,\quad \xi^{(j+1)}:=\xi^\ppu-j\frac{\xi}{3p}\nonumber\\
&&\r_1^\ppj:=\r^\ppj\ ,\quad s_1^\ppj:=s^\ppj-\frac{s}{9p}\ ,\quad \d_1^\ppj:=\d^\ppj-\frac{\d}{9p}\nonumber\\
&&r_1^\ppj:=r^\ppj\ ,\quad  \xi_1^\ppj:=\xi^\ppj\ ,\quad {\cal X}_j:=\sup\{|x|:\ x\in {\mathbb X}_{\xi_j}\}
\end{eqnarray*}
with $1\le j\le p$.

\noindent
We assume that for a certain $1\le i\le p$ and all $1\le j\le i$, we have $\HH_j\in \cO_{\r^\ppj,s^\ppj,\d^\ppj, r^\ppj,\xi^\ppj,}$ of the form
\begin{eqnarray}\label{Hi}
&&\HH_j=\hh_0+g_{j-1}+f_j\ , \quad g_{j-1}\in \cN_{\r^\ppj,s^\ppj,\d^\ppj, r^\ppj,\xi^\ppj,}\ ,\quad g_{j-1}-g_{j-2}=\ovl f_{j-1}\\
&&\label{i+1 ineq} \|f_j\|_{\r^\ppj,s^\ppj,\d^\ppj, r^\ppj,\xi^\ppj,}\le \frac{\|f_1\|_{\r^\ppu,s^\ppu,\d^\ppu, r^\ppu,\xi^\ppu,}}{2^{j-1}}\end{eqnarray}
with  $g_{-1}\equiv 0$,  $g_0=\ovl f$. We  want to prove that, if $i<p$,
Lemma~\ref{iterative lemma} can be applied once again, so as to conjugate $\HH_i$ to a suitable $\HH_{i+1}$ such that~\eqref{Hi}--\eqref{i+1 ineq}
are true with $j=i+1$.
To this end, 
according to the discussion in Remark~\ref{stronger iterative lemma},  
we check the  stronger inequalities
\begin{eqnarray}\label{ineq0}
{2{\r'_i}<\r^\ppi\ ,\quad 2r'_i<r^\ppi\ ,\quad  2\xi'_i<\xi^\ppi}\end{eqnarray}
\begin{eqnarray}\label{smallness0}
&&{\cal X}_i\left\|\Im\frac{\omega }{v }\right\|_{\r_i, r_i}<s'_i%=\frac{s}{6p}
\ ,\quad
{\cal X}_i\left\|\frac{{\omega'} }{v }\right\|_{\r_i, r_i}<
\frac{\d'_i}{\d_i}\\
\label{smallness i}
&& \widetilde{\rm c}\frac{{\cal X}_i}{d_i} % \|\frac{1}{v }\|_{\r_i, r_i}
\left\|\frac{f_i}{v }\right\|_{\r_i,s_i,\d_i, r_i,\xi_i}<1\ .
\end{eqnarray}
where
$d_i:=\min\{ \r_i' s_i', r_i'\xi_i', {\d'}_i^2 \}$.
Conditions~\eqref{ineq0} and~\eqref{smallness0} are certainly  verified, since in fact they are implied by the definitions above 
(using also $\d_i\le \frac{2}{3}\d$, ${\cal X}_i\le {\cal X}$) and the two former inequalities in~\eqref{normal form assumptions}.
 To check the validity of~\eqref{smallness i}, we firstly observe that
\begin{align}\label{di}d_i=\min\{\r'^\ppi s'^\ppi,\ (\d'^\ppi)^2,\ r'^\ppi\xi'_j,\}\ge \frac{\rm d}{{81}p^2}\ .\end{align}
Using then ${\rm c}>{162} \widetilde{\rm c}$,${\cal X}_i<{\cal X}$, Equation~\eqref{f1***}, the inequality in~\eqref{i+1 ineq} with $j=i$ and the last inequality in~\eqref{normal form assumptions}, we easily conclude
\begin{eqnarray*}
&&\|f_i\|_{\r^\ppi,s^\ppi,\d^\ppi, r^\ppi,\xi^\ppi}\le \|f_1\|_{\r^\ppu,s^\ppu,\d^\ppu, r^\ppu,\xi^\ppu,}\le {81}\widetilde{\rm c} \frac{{\cal X}}{\rm d}   %\|\frac{1}{v }\|_{\r, r} 
\left\|\frac{\widetilde f}{v }\right\|_{\r, s, \d, r, \xi}\|{ f}\|_{\r, s, \d, r, \xi}
 \nonumber\\
 &&\le
 \frac{1}{\widetilde{\rm c}}\frac{\rm d}{{81} p^2}\frac{1} {{\cal X}}  \left\|\frac{1}{v }\right\|_{\r, r}^{-1}\le  \frac{1}{\widetilde{\rm c}}\frac{d_i} {{\cal X}_i}  \left\|\frac{1}{v }\right\|_{\r^\ppi, r^\ppi}^{-1}\end{eqnarray*}
 which implies~\eqref{smallness i}.

\noindent
 Then the Iterative Lemma is applicable to $\HH_i$, which is conjugated to
 $$\HH_{i+1}=\hh_0+g_{i}+f_{i+1}\ , \quad g_{i}\in \cN_{\r^{(i+1)},s^{(i+1)},\d^{(i+1)}, r^{(i+1)},\xi^{(i+1)}}\ ,\quad g_{i}-g_{i-1}=\ovl f_{i}$$
 with $g_{i}$, $f_{i+1}$ satisfying~\eqref{Hi}  with $j=i+1$.
We prove that~\eqref{i+1 ineq} holds   with $j=i+1$, so as to complete the inductive step. By the thesis of the Iterative Lemma,
\begin{eqnarray*}
\left\|f_{i+1}\right\|_{\r^{(i+1)}, s^{(i+1)}, \d^{(i+1)}, r^{(i+1)}, \xi^{(i+1)}} &\le& \widetilde{\rm c}\Big(\frac{{\cal X}_i}{{d}_i}\left \|\frac{\widetilde f_i}{v }\right\|_{\r^{(i)}, s^{(i)}, \d^{(i)}, r^{(i)}, \xi^{(i)}}\| f_i\|_{\r^{(i)}, s^{(i)}, \d^{(i)}, r^{(i)}, \xi^{(i)}}\nonumber\\
&+& \|\{\phi_i, g_{i-1}\}\|_{\r^\ppi_1-{\r'}^\ppi, s^\ppi_1-{s'}^\ppi, \d^\ppi_1-{\d'}^\ppi, r^\ppi_1-{r'}^\ppi, \xi^\ppi_1-{\xi'}^\ppi}\Big)
\end{eqnarray*}
On the other hand, using~\eqref{f1***},~\eqref{di} and the last assumption in~\eqref{normal form assumptions} and~\eqref{i+1 ineq}    with $j=i$, we obtain
\begin{eqnarray*}
\frac{{\cal X}_i}{{d}_i}\left \|\frac{\widetilde f_i}{v }\right\|_{\r^{(i)}, s^{(i)}, \d^{(i)}, r^{(i)}, \xi^{(i)}}&\le & \frac{81p^2\chi_0}{{\rm d}}\left\|\frac{1}{v }\right\|\left\|\widetilde f_i\right\|\le
\frac{81p^2\chi_0}{{\rm d}}\left\|\frac{1}{v }\right\|\left\| f_i\right\|
\nonumber\\
&\le & \frac{81p^2\chi_0}{{\rm d}}\left\|\frac{1}{v }\right\|\left\| f_1\right\|\nonumber\\
&\le & \frac{(81)^2 \widetilde{\rm c} p^2\chi^2_0}{{\rm d}^2}\left\|\frac{1}{v }\right\|\left\|\frac{\widetilde f}{v }\right\|_{\r, s, \d, r, \xi}\| f\|_{\r, s, \d, r, \xi}
\nonumber\\
&\le&\su{6\widetilde{\rm c}}
\end{eqnarray*}
Furthermore, using~\eqref{geometric series} with $j=1$,~\eqref{bound on phi} and
\begin{eqnarray*}
\|g_0\|_{\r^\ppi_1-{\r'}^\ppi, s^\ppi_1-{s'}^\ppi, \d^\ppi_1-{\d'}^\ppi, r^\ppi_1-{r'}^\ppi, \xi^\ppi_1-{\xi'}^\ppi}&\le& \|f\|_{\r^\ppi_1-{\r'}^\ppi, s^\ppi_1-{s'}^\ppi, \d^\ppi_1-{\d'}^\ppi, r^\ppi_1-{r'}^\ppi, \xi^\ppi_1-{\xi'}^\ppi} \nonumber\\
&\le&\|f\|_{\r^\ppo, s^\ppo, \d^\ppo, r^\ppo, \xi^\ppo}
\end{eqnarray*}
and, with $j=1$, $\cdots$, $i-1$
\begin{eqnarray}\label{deltag}
\|g_j-g_{j-1}\|_{\r^\ppi_1-{\r'}^\ppi, s^\ppi_1-{s'}^\ppi, \d^\ppi_1-{\d'}^\ppi, r^\ppi_1-{r'}^\ppi, \xi^\ppi_1-{\xi'}^\ppi}&\le& \|f_j\|_{\r^\ppi_1-{\r'}^\ppi, s^\ppi_1-{s'}^\ppi, \d^\ppi_1-{\d'}^\ppi, r^\ppi_1-{r'}^\ppi, \xi^\ppi_1-{\xi'}^\ppi} \nonumber\\
&\le&\|f_j\|_{\r^\ppj,s^\ppj,\d^\ppj, r^\ppj,\xi^\ppj,}\nonumber\\
&\le&\frac{\|f_1\|_{\r^\ppu,s^\ppu,\d^\ppu, r^\ppu,\xi^\ppu,}}{2^{j-1}}
\end{eqnarray}
we obtain
\begin{eqnarray*}
&&\|\{\phi_i, g_{i-1}\}\|_{\r^\ppi_1-{\r'}^\ppi, s^\ppi_1-{s'}^\ppi, \d^\ppi_1-{\d'}^\ppi, r^\ppi_1-{r'}^\ppi, \xi^\ppi_1-{\xi'}^\ppi}\nonumber\\
&&\qquad\le \|\{\phi_i, g_{0}\}\|_{\r^\ppi_1-{\r'}^\ppi, s^\ppi_1-{s'}^\ppi, \d^\ppi_1-{\d'}^\ppi, r^\ppi_1-{r'}^\ppi, \xi^\ppi_1-{\xi'}^\ppi}\nonumber\\
&&\qquad\qquad +\sum_{j=1}^{i-1}\|\{\phi_i, g_{j}-g_{j-1}\}\|_{\r^\ppi_1-{\r'}^\ppi, s^\ppi_1-{s'}^\ppi, \d^\ppi_1-{\d'}^\ppi, r^\ppi_1-{r'}^\ppi, \xi^\ppi_1-{\xi'}^\ppi}\nonumber\\
&&\qquad\qquad\le 2\frac{\ovl{\rm c}\chi_i}{d_i}\left\|\frac{\widetilde f_i}{v }\right\|_{\r^\ppi, s^\ppi, \d^\ppi, r^\ppi, \xi^\ppi}\Big(\frac{1}{p}
\|g_0\|_{\r^\ppi_1-{\r'}^\ppi, s^\ppi_1-{s'}^\ppi, \d^\ppi_1-{\d'}^\ppi, r^\ppi_1-{r'}^\ppi, \xi^\ppi_1-{\xi'}^\ppi}\nonumber\\
&&\qquad\qquad +\sum_{j=1}^{i-1}\|g_{j}-g_{j-1}\|_{\r^\ppi_1-{\r'}^\ppi, s^\ppi_1-{s'}^\ppi, \d^\ppi_1-{\d'}^\ppi, r^\ppi_1-{r'}^\ppi, \xi^\ppi_1-{\xi'}^\ppi}
\Big)\nonumber\\
&&\qquad\qquad\le 162 p^2\frac{\ovl{\rm c}\chi_i}{{\rm d}}\left\|\frac{1}{v }\right\|%_{\r^\ppi, s^\ppi, \d^\ppi, r^\ppi, \xi^\ppi}
\left(\frac{\|f\|_{\r, s, \d, r, \xi}}{p}+2\times {81}\widetilde{\rm c} \frac{{\cal X}}{\rm d}   %\|\frac{1}{v }\|_{\r, r} 
\left\|\frac{\widetilde f}{v }\right\|_{\r, s, \d, r, \xi}\|{ f}\|_{\r, s, \d, r, \xi}
\right)\nonumber\\
&&\qquad\qquad\times
\left\|f_i\right\|_{\r^\ppi, s^\ppi, \d^\ppi, r^\ppi, \xi^\ppi}\le \su3\left\|f_i\right\|_{\r^\ppi, s^\ppi, \d^\ppi, r^\ppi, \xi^\ppi}
\end{eqnarray*}
where the last summand is to taken into account only when $i\ge 2$ and we have used that, when $j=0$, $d_i$ can be replaced with $p d_i$. Collecting the bounds above, we find~\eqref{i+1 ineq} holds   with $j=i+1$.
Then the Iterative Lemma can be applied $p$ times, and we get~\eqref{thesis} as a direct consequence of~\eqref{f1***} and~\eqref{deltag} and of the validity of~\eqref{i+1 ineq} for all $1\le i\le p$. %The inequalities in~\eqref{phi close to id***}
%  follow from~\eqref{phi close to id} and~\eqref{values0}. For example,
%  \begin{align}
%  | I _*- I _0|\le\sum_{j=1}^p| I ^{(j)}- I ^{(j-1)}|\le2\frac{\chi_{0}}{{s'}^{(0)}}\left\|\frac{f_0}{v }\right\|+ 4\frac{\chi_{1}}{{s'}^{(1)}}\left\|\frac{1}{v }\right\|\|f_1\|\le 19\frac{\chi_{0}}{{s}^{(0)}}\left\|\frac{f_0}{v }\right\|
%  \end{align}
%etc. 
$\qquad \square$

\section{Theorem~B}\label{Theorem B}
Let us consider the vector--field~\eqref{perturbed}.

\subsection*{Functional Setting and Statement}
For ease of notation, we discuss the case $n=m=1$. The generalization is left as an exercise.\\
Let us then consider a 3--dimensional vector--field
 \begin{eqnarray*} ({\rm I}, y, \varphi )\in{\mathbb P}_{r, \sigma, s}:={\mathbb I}_r \times{\mathbb Y}_\sigma \times {{\mathbb T}}_s\to X=(X_1, X_2, X_3)\in{\mathbb C}^3\end{eqnarray*} 
 where ${\mathbb I}\subset {\mathbb R}$, ${\mathbb Y}\subset {\mathbb R}$ are open and connected; ${\mathbb T}= {\mathbb R}/(2\pi {\mathbb Z})$,
 which has the form~\eqref{perturbed}.

  \noindent
We assume each $X_i$ to be holomorphic in
${\mathbb P}_{r, \sigma, s}$, meaning the it has a finite weighted norm defined below. If this holds, we simply write $X\in {\cal O}^3_{r, \sigma, s}$.

\noindent
For functions $f:\ ({\rm I}, y, \varphi )\in{\mathbb I}_r\times {\mathbb Y}_\sigma\times {\mathbb T}_s\to {\mathbb C}$, we write $f\in{\cal O}_{r, \sigma, s}$ if $f$ is holomorphic in ${\mathbb P}_{r, \sigma, s}$. We let
\begin{eqnarray*}\|f\|_{u}:=\sum_{k\in {\mathbb Z}}\,\sup_{{\mathbb I}_r\times{\mathbb Y}_\sigma}
|f_{\kappa}({\rm I}, y)|
\,e^{|k|s}\qquad u=(r, \sigma, s)
\end{eqnarray*}
where
\begin{eqnarray*}f=\sum_{k\in {\mathbb Z}} f_k({\rm I}, y)e^{{\rm i} k\varphi }\end{eqnarray*}
is the Fourier series associated to $f$ relatively to the $\varphi $--coordinate.  For $\varphi $--independent 
  functions or vector--fields we simply write $\|\,\cdot\,\|_{r, \sigma}$.
\\
For vector--fields $X:\ ({\rm I}, y, \varphi )\in{\mathbb I}_r\times {\mathbb Y}_\sigma\times {\mathbb T}_s\to X=(X_1, X_2, X_3)\in{\mathbb C}^3$, we write $X\in{\cal O}^3_{r, \sigma, s}$ if $X_i\in{\cal O}_{r, \sigma, s}$ for $i=1$, $2$, $3$. We 
define the {\it weighted norms} 
\begin{eqnarray*}\VERT X \VERT_{u}^{w}:=\sum_i w^{-1}_i\|X_i\|_{u}\end{eqnarray*}
where
$w=(w_1$, $w_2$, $w_3)\in {\mathbb R}_+^3$ are the {\it weights}.

\vskip.1in
\noindent
We shall prove the following

\vskip.1in
\noindent
{\bf Theorem~B} {\it
Let $u=(r, \sigma, s)$; $X=N+P\in {\cal O}^3_{u}$ and let $w=(\rho$, $\tau$, $t)\in {\mathbb R}_+^3$.
Put
\begin{eqnarray*}
 Q:=3\,{\rm diam}({\mathbb Y}_\sigma)\left\|\frac{1}{v}\right\|_{r, \sigma}\end{eqnarray*}
 and\footnote{${\rm diam}({\cal A})$ denotes diameter of the set ${\cal A}$.} assume that for some ${p}\in {\mathbb N}$,  $s_2\in {\mathbb R}_+$, the following inequalities are satisfied:
 \begin{align}\label{NEWu+positive}
0<\rho<\frac{r}{8}\, ,\quad 0<\tau< e^{-s_2}\frac{\sigma}{8}\, ,\quad 0<t<\frac{s}{10}
\end{align}
and
 \begin{align}\label{NEWnewcond2}\chi&:= \frac{{\rm diam}({\mathbb Y}_\sigma)}{s_2}\left\|\frac{\partial_y v}{v}\right\|_{r, \sigma}
\le 1	\\\label{theta1}\theta_1&:= 2\,e^{s_2}{\rm diam}({\mathbb Y}_\sigma)\left\|\frac{\partial_y\omega}{v}\right\|_{r, \sigma}\frac{\tau}{t}\le 1\\
\theta_2&:= 4\,{\rm diam}({\mathbb Y}_\sigma)\left\|\frac{\partial_{\rm I} v}{v}\right\|_{r, \sigma}\frac{\rho}{\tau}\le 1\\
\theta_3&:= 8\,{\rm diam}({\mathbb Y}_\sigma)
\left\|\frac{\partial_{\rm I}\omega}{v}\right\|_{r, \sigma}\frac{\rho}{t}\le 1\\
\label{NEWnewsmallness}\eta^2&:= \max\left\{\frac{{\rm diam}({\mathbb Y}_\sigma)}{t}\left\|\frac{\omega}{v}\right\|_{r, \sigma}\, ,\ 2^7\,e^{2 s_2}Q^2  (\VERT P\VERT_{u}^{w})^2\right\}<\frac{1}{{p}}\, . \end{align}
Then, with
$$u_*=(r_\star, \sigma_\star, s_\star)\, ,\quad r_\star:=r-8\rho\, , \quad\sigma_\star=\sigma-8 e^{s_2}\tau\, , \quad s_\star=s-10 t$$
 there exists  a real--analytic change of coordinates $\Phi_\star$
such that $X_\star:=\Phi_\star X\in {\cal O}^3_{u_\star}$ and
 $X_\star=N+P_\star$,
with
$$\VERT P_\star\VERT^w_{u_\star}<2^{-({p}+1)}\VERT P\VERT^w_{u}\, . $$
}

\subsection*{The Lie Series}

  \begin{definition}\label{Vector--field time--one  flow}\rm For a given $C^{\infty}$ vector--field $Y$, we denote as
$${\cal L}_Y:=[Y, \cdot]$$
 the {\it vector--field Lie operator} (or, simply, {\it Lie operator}, when there is not risk of confusion), where $$[Y, X]:=J_X Y-J_Y X\, ,\quad (J_X)_{ij}:=\partial_{x_j} X_i$$denotes the {\it vector--field Lie brackets} (respectively, {\it Lie brackets}) of two vector--fields.
The map \begin{align}\label{Lie}e^{{\cal L}_Y}:=\sum_{{{k}}=0}^{+\infty}\frac{{\cal L}^{{k}}_Y}{{{k}}!}\end{align}is called {\it vector--field Lie series} (respectively, {\it Lie series}) generated by $Y$.
 \end{definition}

 \vskip.1in
 \noindent
 The importance of  the Lie operator and the Lie series relies in the following

\begin{proposition}\label{prop: Lie operator}
Assume that %, for all $\tau$ in a neighbourhood of $\tau=0$,  $y\to \Phi^Y_\t(y)$  is a diffeomorphism of some open set of ${\mathbb R^n}$ and that
% the Lie operator $e^{{\cal L}_Y} $ converges uniformly. Then 
$e^{{\cal L}_Y}$  is well defined. Then the time--one map of $Y$,
 $\Phi^Y_1$, carries the ODE~\eqref{oldVF} to~\eqref{newVF}, where
$Z=e^{{\cal L}_Y} X$.
\end{proposition}
Proposition~\ref{prop: Lie operator} is well--known. It is however not so easy to find a proof 
not based on differential geometry tools. For this reason, we defer a direct proof of it to Appendix~\ref{app: Lie operator}.

\subsection*{Convergence of the Lie Series}
The purpose of this section is to prove the following
\begin{proposition}\label{Lie Series}
Let $0<w<u$, $Y\in {\cal O}^3_{u+w}$, 
\begin{align}\label{q}q:=3\VERT Y\VERT^w_{u+w}<1\, . \end{align}
Then the Lie series $e^{{\cal L}_Y}$
defines an operator
$$e^{{\cal L}_Y}:\quad  {\cal O}^3_{u}\to {\cal O}^3_{u-w}$$
and its tails  $$\Phi_h:=\sum_{{{k}}= h}^{+\infty}\frac{{\cal L}^{{k}}_Y}{{{k}}!}$$ verify
$$\left\VERT e^{{\cal L}_Y}_m W\right\VERT^w_{u-w}\le \frac{q^m}{1-q}\VERT W\VERT_{u}^w\qquad \forall\ W\in {\cal O}^3_{u}\, . $$
\end{proposition}

\begin{lemma}
The wighted norm affords the following obvious properties.

\begin{itemize}
\item[{\tiny\textbullet}] {Monotonicity:} 
\begin{align}\label{monotonicity}
\VERT X \VERT_{u}^{w}\le \VERT X \VERT_{u'}^{w}\, ,\quad \VERT X \VERT_{u}^{w'}\le \VERT X \VERT_{u}^{w}\quad \forall\ u\le u'\, ,\ w\le w'
\end{align}
where $u\le u'$ means $u_i\le u_i'$ for $i=1$, $2$, $3$.

\item[{\tiny\textbullet}] {Homogeneity:} 
\begin{align}\label{homogeneity}\VERT X \VERT_{u}^{\alpha w}=\alpha^{-1}\VERT X \VERT_{u}^{w}\qquad \forall\ \alpha>0\, . \end{align}
\end{itemize}
\end{lemma}

\begin{lemma}\label{Lie brackets}
Let  $w<u\le u_0$;
$Y\in {\cal O}^3_{u_0}$, $W\in {\cal O}^3_{u}$. Then
$$\VERT{\cal L}_Y[W]\VERT^{u_0-u+w}_{u-w}\le \VERT Y\VERT^{w}_{u-w}\VERT W\VERT ^{u_0-u+w}_{u}+\VERT W\VERT^{u_0-u+w}_{u-w}\VERT Y\VERT^{u_0-u+w}_{u_0}\, . $$
\end{lemma}
{\bf Proof\ }
One has 
\begin{eqnarray*}
\VERT{\cal L}_Y[W]\VERT^{u_0-u+w}_{u-w}&=&\VERT J_W Y-J_Y W\VERT^{u_0-u+w}_{u-w}\nonumber\\
&\le&  \VERT J_W Y\VERT^{u_0-u+w}_{u-w}+\VERT J_Y W\VERT^{u_0-u+w}_{u-w}
\end{eqnarray*}
Now,  $(J_W Y)_i=\partial_{{\rm I}} W_i Y_1+\partial_{y} W_i Y_2+\partial_{\varphi } W_i Y_3$, so, using Cauchy inequalities,
\begin{eqnarray*}
 \|(J_W Y)_i\|_{u-w}&\le&  \|\partial_{{\rm I}} W_i\|_{u-w}\| Y_1\|_{u-w}+\|\partial_{y} W_i\|_{u-w} \|Y_2\|_{u-w}+\|\partial_{\varphi } W_i\|_{u-w} \|Y_3\|_{u-w}\nonumber\\
 &\le& w_1^{-1}\| W_i\|_{u}\| Y_1\|_{u-w}+w_2^{-1}\|W_i\|_{u} \|Y_2\|_{u-w}+w_3^{-1}\|W_i\|_{u} \|Y_3\|_{u-w}\nonumber\\
 &=&\VERT Y\VERT^w_{u-w}\|W_i\|_{u} 
\end{eqnarray*}
Similarly,
\begin{eqnarray*}\|(J_Y W)_i\|_{u-w}\le \VERT W\VERT_{u-w}^{u_0-u+w}\|Y_i\|_{u_0}\, .  \end{eqnarray*}
Taking the $u_0-u+w$--weighted norms, the thesis follows. $\quad \square$
\begin{lemma}\label{iterateL}
Let $0<w<u\in {\mathbb R}^3$, $Y\in {\cal O}^3_{u+w}$, $W\in {\cal O}^3_{u}$. Then
$$\VERT{\cal L}^{{k}}_Y[W]\VERT^{w}_{u-w}\le 3^{{k}} {{k}}!\left(\VERT Y\VERT ^w_{u+w}\right)^{{k}}\VERT W\VERT^w _{u-w}\, . $$
\end{lemma}

\proof\, We apply Lemma~\ref{Lie brackets} with $W$ replaced by ${\cal L}^{i-1}_Y[W]$, 
$u$ replaced by $u-(i-1)w/{{k}}$, 
$w$ replaced by $w/{{k}}$ and, finally, $u_0=u+w$. With $\VERT\cdot\VERT_i^{w}=\VERT\cdot\VERT^{w}_{u-i\frac{w}{{{k}}}}$, $0\le i\le {{k}}$, so that $\VERT\cdot\VERT^w_0=\VERT\cdot\VERT^w_{u}$ and $\VERT\cdot\VERT^w_{{k}}=\VERT\cdot\VERT^w_{u-w}$,
\begin{eqnarray*}
\VERT{\cal L}^i_Y[W]\VERT^{w+w/{{k}}}_i&=&\left\VERT\left[Y, {\cal L}^{i-1}_Y[W]\right]\right\VERT^{w+w/{{k}}}_i\nonumber\\
&\le&  \VERT Y\VERT^{w/{{k}}}_{i}\VERT{\cal L}^{i-1}_Y[W]\VERT^{w+w/{{k}}}_{i-1}+
\VERT Y\VERT^{w+w/{{k}}}_{u+w}\VERT{\cal L}^{i-1}_Y[W]\VERT^{w+w/{{k}}}_{i}\,.
\end{eqnarray*}
Hence, de--homogenizating,
\begin{eqnarray*}
\frac{{{k}}}{{{k}}+1}\VERT{\cal L}^i_Y[W]\VERT^{w}_i&\le&  {{k}} \frac{{{k}}}{{{k}}+1}\VERT Y\VERT^{w}_{i}\VERT{\cal L}^{i-1}_Y[W]\VERT^{w}_{i-1}+
\frac{{{k}}^2}{({{k}}+1)^2}\VERT Y\VERT^{w}_{u+w}\VERT{\cal L}^{i-1}_Y[W]\VERT^{w}_{i}\nonumber\\
&\le& \frac{{{k}}^2}{{{k}}+1}\left(1+\frac{1}{{{k}}+1}\right)\VERT Y\VERT^{w}_{u+w}\VERT{\cal L}^{i-1}_Y[W]\VERT^{w}_{i-1}
\end{eqnarray*}
Eliminating the common factor $\frac{{{k}}}{{{k}}+1}$ and 
 iterating ${{k}}$ times from $i={{k}}$, by Stirling, we get
\begin{eqnarray*}
\VERT {\cal L}^{{k}}_Y[W]\VERT^w _{u-w}
&\le&   {{k}}^{{k}}\left(1+\frac{1}{{k}}\right)^{{k}}\left(\VERT Y\VERT ^w_{u+w}\right)^{{k}}\VERT W\VERT^w _{u-w}\nonumber\\
&\le&  e^{{k}} {{k}}!\left(\VERT Y\VERT ^w_{u+w}\right)^{{k}}\VERT W\VERT^w _{u-w}\nonumber\\
&< &3^{{k}} {{k}}!\left(\VERT Y\VERT ^w_{u+w}\right)^{{k}}\VERT W\VERT^w _{u-w}
\end{eqnarray*}
as claimed. $\quad \square$

\vskip.1in
\noindent
{\bf Proof of Proposition~\ref{Lie Series}} Immediate from Lemma~\ref{iterateL}. $\qquad \square$

\subsection*{The Homological Equation}
\begin{proposition}\label{homeq1}
Let 
\begin{align}\label{NZ}N=(0, v({\rm I}, y), \omega({\rm I}, y))\, ,\qquad Z=(Z_1({\rm I}, y, \varphi ), Z_2({\rm I}, y, \varphi ), Z_3({\rm I}, y, \varphi ))\end{align}
belong to ${\cal O}^3_{r, \sigma, s}$ and assume~\eqref{existence}.
Then  the ``homological equation''
\begin{align}\label{homeq}{\cal L}_N[Y]=Z\end{align}
has a solution $Y\in {\cal O}_{r, \sigma, s-3s_1}$ verifying
\begin{align}\label{bounds}
\VERT Y\VERT_{r, \sigma, s-3 s_1}^{\rho_*, \tau_*, t_*}\le {\rm diam}({\mathbb Y}_\sigma)\left\|\frac{1}{v}\right\|_{r, \sigma}\VERT Z\VERT_{r, \sigma, s}^{\rho, \tau, t}
\end{align}
with $\rho_*$, $\tau_*$, $t_*$ as in~\eqref{bounds1}.
\end{proposition}

\noindent
Fix $y_0\in {\mathbb Y}$; $v$, $\omega:{\mathbb I}\times {\mathbb Y}\to {\mathbb R}$, with $v\not\equiv  0$. We define, formally, the operators ${\cal F}_{v, \omega}$ and ${\cal G}_{v, \omega}$
as acting on  functions $g:{\mathbb I}\times {\mathbb Y}\times {\mathbb T}\to {\mathbb R}$ as
\begin{align}\label{FandG}
{\cal F}_{v, \omega}[g]({\rm I}, y, \varphi )&:= \int_{y_0}^y\frac{g\left({\rm I}, \eta, \varphi +\int_y^\eta \frac{\omega({\rm I}, \eta')}{v({\rm I}, \eta')}d\eta' \right)}{v({\rm I}, \eta)}d\eta\nonumber\\
{\cal G}_{v, \omega}[g]({\rm I}, y, \varphi )&:= \int_{y_0}^y\frac{g\left({\rm I}, \eta, \varphi +\int_y^\eta \frac{\omega({\rm I}, \eta')}{v({\rm I}, \eta')}d\eta' \right)e^{-\int_y^\eta \frac{\partial_y v({\rm I}, \eta')}{v({\rm I}, \eta')}d\eta'}}{v({\rm I}, \eta)}d\eta
\end{align}

\noindent
Observe that, when existing,
${\cal F}_{v, \omega}$, ${\cal G}_{v, \omega}$  send zero--average functions to zero--average functions.

\noindent
The existence ${\cal F}_{v, \omega}$, ${\cal G}_{v, \omega}$ is established by the following 

\begin{lemma}\label{estimates}
If inequalities~\eqref{existence} hold,
then
$${\cal F}_{v, \omega}\, ,\ {\cal G}_{v, \omega}:\quad {\cal O}_{r, \sigma, s}\to {\cal O}_{r, \sigma, s-s_1}$$
and
$$\|{\cal F}_{v, \omega}[g]\|_{r, \sigma, s-s_1}\le {\rm diam}({\mathbb Y}_\sigma)\left\|\frac{g}{v}\right\|_{r, \sigma, s}\, ,\quad \|{\cal G}_{v, \omega}[g]\|_{r, \sigma, s-s_1}\le e^{s_2}\,{\rm diam}({\mathbb Y}_\sigma) \left\|\frac{g}{v}\right\|_{r, \sigma, s}$$
 \end{lemma}
 The proof of Lemma~\ref{estimates} is obvious from the definitions~\eqref{FandG}.

\proof\,{\bf of Proposition~\ref{homeq1}\ }We expand $Y_j$ and $Z_j$ along the Fourier basis
$$Y_j({\rm I}, y, \varphi )=\sum_{k\in {\mathbb Z}} Y_{j, k}({\rm I}, y)e^{{\rm i} k\varphi }\, ,\quad Z_j({\rm I}, y, \varphi )=\sum_{k\in {\mathbb Z}} Z_{j, k}({\rm I}, y)e^{{\rm i} k\varphi }\, ,\quad  j=1,\ 2,\ 3$$
Using $${\cal L}_N[Y]=[N, Y]=J_Y N-J_N Y$$
where
$(J_Z)_{ij}=\partial_j Z_i$
are the Jacobian matrices,
we rewrite~\eqref{homeq} as

\begin{align}\label{equations}
Z_{1, k}({\rm I}, y)&=v({\rm I}, y)\partial_y Y_{1, k} +{\rm i} k \omega({\rm I}, y) Y_{1, k} \nonumber\\
Z_{2, k}({\rm I}, y)&=v({\rm I}, y)\partial_y Y_{2, k} +({\rm i} k \omega({\rm I}, y)-\partial_y v({\rm I}, y)) Y_{2, k}-\partial_{\rm I} v({\rm I}, y) Y_{1, k}\nonumber\\
Z_{3, k}({\rm I}, y)&=v({\rm I}, y)\partial_y Y_{3, k} +{\rm i} k \omega({\rm I}, y) Y_{3, k}-\partial_{\rm I} \omega({\rm I}, y) Y_{1, k}-\partial_y \omega({\rm I}, y) Y_{2, k}\,.
\end{align}
Regarding~\eqref{equations} as equations for  $Y_{j, k}$, we find the solutions
\begin{eqnarray*}
Y_{1, k}&=&\int_{y_0}^y\frac{Z_{1, k}({\rm I}, \eta)}{v({\rm I}, \eta)}e^{{\rm i} k\int_y^\eta \frac{\omega({\rm I}, \eta')}{v({\rm I}, \eta')}d\eta'}d\eta\nonumber\\
Y_{2, k}&=&\int_{y_0}^y\frac{Z_{2, k}({\rm I}, \eta)+\partial_{\rm I} v
Y_{1, k}
}{v({\rm I}, \eta)}e^{\int_y^\eta \frac{{\rm i} k \omega({\rm I}, \eta')-\partial_y v({\rm I}, \eta')}{v({\rm I}, \eta')}d\eta'}d\eta\nonumber\\
Y_{3, k}&=&\int_{y_0}^y\frac{Z_{3, k}({\rm I}, \eta)+\partial_{\rm I} \omega({\rm I}, \eta) Y_{1, k}+\partial_y \omega({\rm I}, \eta) Y_{2, k}}{v({\rm I}, \eta)}e^{{\rm i} k\int_y^\eta \frac{\omega({\rm I}, \eta')}{v({\rm I}, \eta')}d\eta'}d\eta\nonumber\\
\end{eqnarray*} 
multiplying by $e^{{\rm i} k\varphi }$ and summing over $k\in {\mathbb Z}$ we find \begin{eqnarray*}
Y_1&=&{\cal F}_{v, \omega}[Z_1]\nonumber\\
 Y_2&=&{\cal G}_{v, \omega}[Z_2]+{\cal G}_{v, \omega}[\partial_{\rm I} v\,Y_1]\, ,\nonumber\\ 
 Y_3&=&{\cal F}_{v, \omega}[Z_3]+{\cal F}_{v, \omega}[\partial_{\rm I} \omega\,Y_1]+{\cal F}_{v, \omega}[\partial_y \omega\,Y_2]\, . \end{eqnarray*} 
Then, by Lemma~\ref{estimates},
\begin{eqnarray*}
&&\|Y_1\|_{r, \sigma, s-s_1}\le  {\rm diam}({\mathbb Y}_\sigma)\left\|\frac{1}{v}\right\|_{r, \sigma}\left\|Z_1\right\|_{r, \sigma, s}\nonumber\\
&&\|Y_2\|_{r, \sigma, s-2s_1}%&\le  e^{s_2}{\rm diam}({\mathbb Y}_\sigma)\left\|\frac{1}{v}\right\|_{r, \sigma}\left\|Z_2\right\|_{r, \sigma, s-s_1}+e^{s_2}{\rm diam}({\mathbb Y}_\sigma)\left\|\frac{{\partial_{\rm I} v}}{v}\right\|_{r, \sigma}\left\|Y_1\right\|_{r, \sigma, s-s_1}\nonumber\\
\le e^{s_2}{\rm diam}({\mathbb Y}_\sigma)\left\|\frac{1}{v}\right\|_{r, \sigma}\left\|Z_2\right\|_{r, \sigma, s-s_1}+e^{s_2}{\rm diam}({\mathbb Y}_\sigma)^2\left\|\frac{1}{v}\right\|_{r, \sigma}\left\|\frac{\partial_{\rm I} v}{v}\right\|_{r, \sigma}\left\|Z_1\right\|_{r, \sigma, s}\nonumber\\
&&\|Y_3\|_{r, \sigma, s-3s_1}%&\le  {\rm diam}({\mathbb Y}_\sigma)\left\|\frac{1}{v}\right\|_{r, \sigma}\left\|Z_3\right\|_{r, \sigma, s-2s_1}+{\rm diam}({\mathbb Y}_\sigma)\left\|\frac{{\partial_{\rm I} \omega}}{v}\right\|_{r, \sigma}\left\|Y_1\right\|_{r, \sigma, s-2s_1}\nonumber\\
%&+ {\rm diam}({\mathbb Y}_\sigma)\left\|\frac{\partial_y \omega}{v}\right\|_{r, \sigma}\left\|Y_2\right\|_{r, \sigma, s-2s_1}
%\nonumber\\
%&\le  {\rm diam}({\mathbb Y}_\sigma)\left\|\frac{1}{v}\right\|_{r, \sigma}\left\|Z_3\right\|_{r, \sigma, s-2s_1}+{\rm diam}({\mathbb Y}_\sigma)^2\left\|\frac{1}{v}\right\|_{r, \sigma}
%\left\|\frac{{\partial_{\rm I} \omega}}{v}\right\|_{r, \sigma}
%\left\|Z_1\right\|_{r, \sigma, s-s_1}\nonumber\\
%
%&+ {\rm diam}({\mathbb Y}_\sigma)\left\|\frac{{\partial_y \omega}}{v}\right\|_{r, \sigma}\left(e^{s_2}{\rm diam}({\mathbb Y}_\sigma)\left\|\frac{1}{v}\right\|_{r, \sigma}\left\|Z_2\right\|_{r, \sigma, s-s_1}\right.\nonumber\\
%&+ \left.e^{s_2}{\rm diam}({\mathbb Y}_\sigma)^2\left\|\frac{1}{v}\right\|_{r, \sigma}
%\left\|\frac{{\partial_{\rm I} v}}{v}\right\|_{r, \sigma}
%\left\|Z_1\right\|_{r, \sigma, s}\right)\nonumber\\
%
\le  {\rm diam}({\mathbb Y}_\sigma)\left\|\frac{1}{v}\right\|_{r, \sigma}\left\|Z_3\right\|_{r, \sigma, s-2s_1}
+e^{s_2}{\rm diam}({\mathbb Y}_\sigma)^2\left\|\frac{1}{v}\right\|_{r, \sigma}\left\|\frac{\partial_y\omega}{v}\right\|_{r, \sigma}
\left\|Z_2\right\|_{r, \sigma, s-s_1}
\nonumber\\
&&\qquad+{\rm diam}({\mathbb Y}_\sigma)^2\left\|\frac{1}{v}\right\|_{r, \sigma}
\left(\left\|\frac{\partial_{\rm I}\omega}{v}\right\|_{r, \sigma}+e^{s_2}{\rm diam}({\mathbb Y}_\sigma)\left\|\frac{\partial_{\rm I} v}{v}\right\|_{r, \sigma} \left\|\frac{\partial_y\omega}{v}\right\|_{r, \sigma}\right)\left\|Z_1\right\|_{r, \sigma, s}
\end{eqnarray*}
Multiplying the inequalities above by $\rho_*^{-1}$, $\tau_*^{-1}$, $t_*^{-1}$ respectively and taking the sum, we find~\eqref{bounds}, with 
\begin{eqnarray*}
\frac{1}{\rho}&=&\frac{1}{\rho_*}+e^{s_2}{\rm diam}({\mathbb Y}_\sigma)\left\|\frac{\partial_{\rm I} v}{v}\right\|_{r, \sigma}\frac{1}{\tau_*}+{\rm diam}({\mathbb Y}_\sigma)
\left(\left\|\frac{\partial_{\rm I}\omega}{v}\right\|_{r, \sigma}+e^{s_2}{\rm diam}({\mathbb Y}_\sigma)\left\|\frac{\partial_{\rm I} v}{v}\right\|_{r, \sigma} \left\|\frac{\partial_y\omega}{v}\right\|_{r, \sigma}\right)\frac{1}{t_*}\nonumber\\
\frac{1}{\tau}&=&\frac{e^{s_2}}{\tau_*}+e^{s_2}{\rm diam}({\mathbb Y}_\sigma)\left\|\frac{\partial_y\omega}{v}\right\|_{r, \sigma}\frac{1}{t_*}\nonumber\\
\frac{1}{t}&=&\frac{1}{t_*}\, . 
\end{eqnarray*}
We recognize that, under conditions~\eqref{w*positive}, $\rho_*$, $\tau_*$, $t_*$ in~\eqref{bounds1} solve the equations above. 
$\qquad \square$

\subsection*{The Step Lemma}
The  Step Lemma is now as follows.
\begin{lemma}\label{iteration lemma NEW}
Let $X=N
+P\in {\cal O}^3_{u}$, with $u=(r, \sigma, s)$, $N$ as in~\eqref{NZ}, %$s_1$, 
$s_2>0$. Assume 
\begin{align}\label{newcond1}\begin{split}
&2\,e^{s_2}{\rm diam}({\mathbb Y}_\sigma)\left\|\frac{\partial_y\omega}{v}\right\|_{r, \sigma}\frac{\tau}{t}\le 1\\
&4\,{\rm diam}({\mathbb Y}_\sigma)\left\|\frac{\partial_{\rm I} v}{v}\right\|_{r, \sigma}\frac{\rho}{\tau}\le 1\\
&8\,{\rm diam}({\mathbb Y}_\sigma)
\left\|\frac{\partial_{\rm I}\omega}{v}\right\|_{r, \sigma}\frac{\rho}{t}\le 1
\end{split}
\end{align}
\begin{align}\label{newcond2}
\frac{{\rm diam}({\mathbb Y}_\sigma)}
{t}\left\|\frac{\omega}{v}\right\|_{r, \sigma}\le 1\, ,\quad \frac{{\rm diam}({\mathbb Y}_\sigma)}{s_2}\left\|\frac{\partial_y v}{v}\right\|_{r, \sigma}
\le 1
\end{align}
\begin{align}
\label{u+positive}
0<\rho<\frac{r}{4}\, ,\quad 0<\tau<\frac{\sigma}{4}e^{-s_2}\, ,\quad 0<t<\frac{s}{5}\, ,\quad w:=(\rho, \tau, t)
\end{align}
\begin{align}
\label{newsmallness} 
2Q\VERT P\VERT ^w_u<1
\end{align}
Let
\begin{align}u_+:=(r-4\rho, \sigma-4\tau e^{s_2}, s-5t)\, . \end{align}
 Then there exists $w_*\in {\mathbb R}_+^3$, $Y\in {\cal O}^3_{u_++w_*}$  such that
$X_+:=e^{{\cal L}_Y}X\in {\cal O}^3_{u_+}$ and 
$X_+=N+
P_+$, with
\begin{align}\label{finalineq}\VERT P_+\VERT_{u_+}^{w}\le 8 e^{s_2}Q
(\VERT P \VERT_u^w)^2\, . \end{align}
\end{lemma}

\proof\, We prove, more\footnote{The  inequality in~\eqref{newcond2} guarantees that one can take $s_1=t$, while the inequalities in~\eqref{newcond1} and~\eqref{u+positive} imply
\begin{align}\frac{1}{\rho_*}\ge \frac{1}{2\rho}\, ,\quad \frac{1}{\tau_*}\ge \frac{e^{-s_2}}{2\tau}
\end{align}
\begin{align}w_*<2 e^{s_2} w\, ,\qquad u_*\ge u_+>0\, . \end{align}
Then~\eqref{finalineq} is implied by~\eqref{P+}, monotonicity and homogeneity~\eqref{monotonicity}--\eqref{homogeneity}, and the inequality in~\eqref{newsmallness}. }  generally, that, for any
$s_1>0$ and provided that 
\begin{align}\label{existence}
\frac{{\rm diam}({\mathbb Y}_\sigma)}{s_1}\left\|\frac{\omega}{v}\right\|_{r, \sigma}\le 1\, ,\quad \frac{{\rm diam}({\mathbb Y}_\sigma)}{s_2}\left\|\frac{\partial_y v}{v}\right\|_{r, \sigma}
\le 1\end{align} 
$$
Q \VERT P\VERT^{w} _{u}<1\qquad Q:=3{\rm diam}({\mathbb Y}_\sigma)\left\|\frac{1}{v}\right\|_{r, \sigma}\, ,\quad w=(\rho, \tau, t)
$$
defining
\begin{align}\label{w*positive}w_*=(\rho_*, \tau_*, t_*)\in {\mathbb R}_+^3\, ,\qquad u_*=(r-2\rho_*, \sigma-2\tau_*, s-3s_1-2t_*)\in {\mathbb R}_+^3\, . \end{align}
with
\begin{eqnarray}\label{bounds1}
\frac{1}{\rho_*}&=&\frac{1}{\rho}-{\rm diam}({\mathbb Y}_\sigma)\left\|\frac{\partial_{\rm I} v}{v}\right\|_{r, \sigma}\left(\frac{1}{\tau}-e^{s_2}{\rm diam}({\mathbb Y}_\sigma)
\left\|\frac{\partial_y\omega}{v}\right\|_{r, \sigma}\frac{1}{t}
\right)\nonumber\\
&-& {\rm diam}({\mathbb Y}_\sigma)
\left(\left\|\frac{\partial_{\rm I}\omega}{v}\right\|_{r, \sigma}+e^{s_2}{\rm diam}({\mathbb Y}_\sigma)\left\|\frac{\partial_{\rm I} v}{v}\right\|_{r, \sigma} \left\|\frac{\partial_y\omega}{v}\right\|_{r, \sigma}\right)\frac{1}{t}\nonumber\\
\frac{1}{\tau_*}&=&\frac{e^{-s_2}}{\tau}-{\rm diam}({\mathbb Y}_\sigma)\left\|\frac{\partial_y\omega}{v}\right\|_{r, \sigma}\frac{1}{t}\nonumber\\
t_*&=&t
\end{eqnarray}
for any $Y\in {\cal O}^3_{u_*+w_*}$, it is \begin{align}\label{P+}
 \VERT P_+\VERT^{w_*}_{u_*}
\le \frac{2Q
\left(\VERT P\VERT^w_u\right)^2
}{1-Q
\left\VERT P\right\VERT^w_u}
\end{align}

\noindent
We look for $Y$ such that $X_+:=e^{{\cal L}_Y} X$ has the desired properties. 
\begin{eqnarray*}
e^{{\cal L}_Y} X&=&e^{{\cal L}_Y}\left(N
+P\right)=N
+P+{\cal L}_Y N+e_2^{{\cal L}_Y} N
+e_1^{{\cal L}_Y}P\nonumber\\
&=&N
+P-{\cal L}_N Y+P_+
\end{eqnarray*}
with
$P_+=e_2^{{\cal L}_Y} N
+e_1^{{\cal L}_Y}P$.
We 
choose $Y$ so that the homological equation 
$${\cal L}_N Y=P$$
is satisfied. By Proposition~\ref{homeq1}, this equation has a solution $Y\in {\cal O}^3_{r, \sigma, s-3s_1}$ verifying 
$$ q:=3\VERT Y\VERT^{w_*} _{r, \sigma, s-3s_1}\le 
3{\rm diam}({\mathbb Y}_\sigma)\left\|\frac{1}{v}\right\|_{r, \sigma}\VERT P\VERT^w_u=Q \VERT P\VERT^w_u<1\, .$$
By Proposition~\ref{Lie Series}, the Lie series $e^{{\cal L}_Y}$
defines an operator
$$e^{{\cal L}_Y}:\quad W\in {\cal O}_{u_*+w_*}\to {\cal O}_{u_*}$$
and its tails $e^{{\cal L}_Y}_m$
verify
\begin{eqnarray*}
\left\VERT e^{{\cal L}_Y}_m W\right\VERT^{w_*}_{u_*}&\le&  \frac{q^m}{1-q}\VERT W\VERT^{w_*}_{u_*+w_*}\nonumber\\
&\le& \frac{\left(Q\VERT P\VERT^w_u\right)^m}{1-Q\VERT P\VERT^w_u}\VERT W\VERT^{w_*}_{u_*+w_*}
\end{eqnarray*}
for all $W\in {\cal O}^3_{u_*+w_*}$.
In particular, $e^{{\cal L}_Y}$ is well defined on ${\cal O}^3_{u}\subset {\cal O}^3_{u_*+w_*}$, hence $P_+\in {\cal O}^3_{u_*}$. The bounds on $P_+$ are obtained as follows. Using the homological equation, one finds
\begin{eqnarray*}
\VERT e_2^{{\cal L}_Y}N\VERT^{w_*}_{u_*}&=& \left\VERT \sum_{{{k}}=1}^{\infty}\frac{{\cal L}^{{{k}}+1}_Y N}{({{k}}+1)!}\right\VERT^{w_*}_{u_*}\nonumber\\
&\le&
\sum_{{{k}}=1}^{\infty}\frac{1}{({{k}}+1)!}\left\VERT {\cal L}^{{{k}}+1}_Y N\right\VERT^{w_*}_{u_*}\nonumber\\
&=&\sum_{{{k}}=1}^{\infty}\frac{1}{({{k}}+1)!}\left\VERT {\cal L}^{{{k}}}_Y P\right\VERT^{w_*}_{u_*}\nonumber\\
&\le&  \sum_{{{k}}=1}^{\infty}\frac{1}{{{k}}!}\left\VERT {\cal L}^{{{k}}}_Y P\right\VERT^{w_*}_{u_*}\nonumber\\
&\le& \frac{Q
\left(\VERT P\VERT^w_u\right)^2
}{1-Q
\left\VERT P\right\VERT^w_u}
\end{eqnarray*}
The bound
$$
 \VERT e_1^{{\cal L}_Y}P\VERT^{w_*}_{u_*}
\le \frac{Q
\left(\VERT P\VERT^w_u\right)^2
}{1-Q
\left\VERT P\right\VERT^w_u}$$
is even more straightforward. $\quad \square$

\subsection*{Iterations (Proof of Theorem~B)}\label{Proof of NFL}
The proof of Theorem~B is obtained -- following~\cite{poschel93} -- via iterate applications of the Step Lemma. 
At the base step, we let\footnote{{With slight abuse of notations, here and during the proof of Theorem~B, the sub--fix $j$ will denote   the value of a given quantity at the $j^{\rm th}$ step  of the iteration.}}
$$X=X_0:=N+P_0\, ,\quad w=w_0:=(\rho,\tau, t)\, ,\quad u=u_0:=(r,\sigma, s)$$
with $X_0=N+P_0\in {\cal O}^3_{u_0}$. We let
$$
Q_0:=3\,{\rm diam}({\mathbb Y}_\sigma)\left\|\frac{1}{v}\right\|_{r, \sigma}$$Conditions~\eqref{newcond1}--\eqref{newsmallness}   are implied by the assumptions~\eqref{theta1}--\eqref{NEWnewsmallness}. We then conjugate $X_0$ to $X_1=N+P_1\in {\cal O}^3_{u_1}$, where
$$u_1=(r-4\rho, \sigma-4\tau e^{s_2}, s-5t)=:(r_1, \sigma_1, s_1)\, . $$
Then we have 
\begin{align}\label{base step}\VERT P_1\VERT_{u_1}^{w_0}\le 8 e^{s_2}Q_0  \left(\VERT P_0\VERT_{u_0}^{w_0}\right)^2\le \frac{1}{2} \VERT P_0\VERT_{u_0}^{w_0}\, .  \end{align}
We assume, inductively, that, for some $1\le j\le {{p}}$, 
we have  
\begin{align}\label{induction}X_{j}=N+P_{j}\in {\cal O}^3_{u_j}\, ,\qquad \VERT P_{j}\VERT_{u_j}^{w_0}<2^{-(j-1)}
\VERT P_{1}\VERT _{u_{1}}^{w_0}
\end{align} where \begin{align}\label{uj}u_j=(r_j, \sigma_j, s_j)\end{align} with
$$r_j:=r_1-4 (j-1)\frac{\rho}{ {{p}}}\, ,\quad \sigma_j:=\sigma_1-4 e^{s_2}(j-1)\frac{\tau}{ {{p}}}\, ,\quad s_j:=s_1-5 (j-1)\frac{t}{ {{p}}}\, . $$ The case $j=1$  trivially   reduces to the identity $\VERT P_{1}\VERT _{u_{1}}^{w_0}=\VERT P_{1}\VERT _{u_{1}}^{w_0}$. 
We  aim to apply {Lemma}~\ref{iteration lemma NEW}  with $u=u_j$ as in~\eqref{uj} and $$w=w_1:=\frac{w_0}{{{p}}}\, ,\qquad \forall\ 1\le j\le {{p}}\, . $$
Conditions~\eqref{newcond1},~\eqref{newcond2} and~\eqref{u+positive} are easily seen to be implied by~\eqref{theta1},~\eqref{NEWnewcond2},~\eqref{NEWu+positive} and the first condition in~\eqref{NEWnewsmallness} combined with the inequality ${{p}}\eta^2<1$, implied by the choice of ${{p}}$. We check condition~\eqref{newsmallness}.
By homogeneity,
$$\VERT P_j\VERT_{u_j}^{w_1}={{p}}\VERT P_j\VERT_{u_j}^{w_0}\le {{p}}\VERT P_1\VERT_{u_1}^{w_0}\le 8{{p}} e^{s_2}Q_0  \left(\VERT P_0\VERT_{u_0}^{w_0}\right)^2$$whence, using
$$
Q_j=3\,{\rm diam}({\mathbb Y}_{\sigma_j})\left\|\frac{1}{v}\right\|_{r_j, \sigma_j}\le Q_0$$
we see that condition~\eqref{newsmallness} is met:
$$2 Q_j\VERT P_j\VERT_{u_j}^{w_1}\le 16{{p}} e^{s_2}Q^2_0  \left(\VERT P_0\VERT_{u_0}^{w_0}\right)^2<1\, . $$
Then the Iterative Lemma can be applied and we get $X_{j+1}=N+P_{j+1}\in {\cal O}^3_{u_{j+1}}$, with $$\VERT P_{j+1}\VERT _{u_{j+1}}^{w_1}\le 8 e^{s_2}  Q_j \left(\VERT P_j\VERT _{u_j}^{w_1}\right)^2\le 8 e^{s_2}  Q_0 \left(\VERT P_j\VERT _{u_j}^{w_1}\right)^2\,.$$
Using homogeneity again to the extreme sides of this inequality and combining it with~\eqref{induction},~\eqref{base step} and~\eqref{NEWnewsmallness}, we get
\begin{eqnarray*}\VERT P_{j+1}\VERT _{u_{j+1}}^{w_{0}}&\le&  8 {{p}} e^{s_2}  Q_0 \left(\VERT P_j\VERT _{u_j}^{w_0}\right)^2
\le 8 {{p}} e^{s_2}  Q_0  \VERT P_1\VERT _{u_1}^{w_0}\VERT P_j\VERT _{u_j}^{w_0}\nonumber\\
&\le&
64{{p}}e^{2s_2}  Q_0^2\left(\VERT P_0\VERT _{u_0}^{w_0}\right)^2
\VERT P_j\VERT _{u_j}^{w_0}
\le
\frac{1}{2}\VERT P_j\VERT _{u_j}^{w_0}\nonumber\\
&<&2^{-j}\VERT P_1\VERT _{u_1}^{w_0}\, . 
 \end{eqnarray*}
After ${{p}}$ iterations, 
$$\VERT P_{{{p}}+1}\VERT _{u_{{{p}}+1}}^{w_{0}}<2^{-{{p}}}\VERT P_1\VERT _{u_1}^{w_0}<2^{-({{p}}+1)}\VERT P_0\VERT _{u_0}^{w_0}$$so we can take $X_\star=X_{{{p}}+1}$, $P_\star=P_{{{p}}+1}$,  $u_\star=u_{{{p}}+1}$. $\qquad \square$

\appendix
\section{Appendix}
\subsection*{Proof of Proposition~\ref{prop: time one map}}\label{app: time one flow}
We want to consider transformations $$y=(p, q)\in U\subset {\mathbb R}^{2n}\to z=(r, s)\in V\subset {\mathbb R}^{2n}$$ where $U$, $V$ are open sets,
corresponding to 
the  time--$\tau$ flow of the Hamiltonian $\phi$ with initial datum $y$.
Namely, $z=\Phi^\phi_\t(y)$ is defined as the solution
of
$$\left\{
\begin{array}{lll}
\partial_{\tau} z=J\partial_{z} \phi(z)\\\\
z|_{\tau=0}=y
\end{array}
\right.$$with $J$
the symplectic unit in~\eqref{J}.

\noindent
Consider the $\tau$--expansion of $\Phi^\phi_\tau (y)$ about a given $\tau_0$:
$$\Phi^\phi_\tau (y)=\Phi^\phi_{\tau_0} (y)+J\partial \phi\big(\Phi^\phi_{\tau_0} (y)\big)(\tau-\tau_0)+o(\tau-\tau_0)\,.$$
Then it is
\begin{eqnarray*}\HH\big(\Phi^\phi_\tau (y)\big)&=&\HH\big(\Phi^\phi_{\tau_0} (y)\big)+\partial\HH\cdot J\partial_z \phi\big(\Phi^\phi_{\tau_0} (y)\big)(\tau-\tau_0)+o(\tau-\tau_0) \nonumber\\
&=&\HH\big(\Phi^\phi_{\tau_0} (y)\big)+{\cal L}_\phi\HH\big(\Phi^\phi_{\tau_0} (y)\big)(\tau-\tau_0)+o(\tau-\tau_0)
\end{eqnarray*}
This shows
$$\frac{d}{d\tau}\HH\big(\Phi^\phi_\tau (y)\big)={\cal L}_\phi\HH\big(\Phi^\phi_{\tau} (y)$$
whence, iterating,
$$\frac{d^k}{d\tau^k}\HH\big(\Phi^\phi_\tau (y)\big)={\cal L}^k_\phi\HH\big(\Phi^\phi_{\tau} (y)\,.$$
This gives, finally,
$$\HH_1(y)=\HH\big(\Phi^\phi_1 (y)\big)=\sum_{k=0}^\infty\frac{1}{k!}\frac{d^k}{d\tau^k}\HH\big(\Phi^\phi_\tau (y)\big)\Big|_{\tau=0}=\sum_{k=0}^\infty\frac{{\cal L}^k_\phi\HH (y)}{k!}=e^{{\cal L}_\phi}\HH\,.\qquad \square$$

\subsection*{Proof of Proposition~\ref{prop: Lie operator}}\label{app: Lie operator}

 In general, a diffeomorphism
$$\Phi:\quad z\in V\to x=\Phi(z)\in U$$ transforms the Equation~\eqref{oldVF} to~\eqref{newVF}
where
\begin{align}\label{new vectorfield}Z(z)=J(z)^{-1}X\big(\Phi (z)\big)\end{align}
with $J(z)$ being the Jacobian matrix of the transformation, i.e., $$J(z)_{hk}=\partial_{z_k}\Phi_h(z)\,,\quad {\rm if}\quad \Phi=(\Phi_1, \ldots, \Phi_n)\,.$$
We want to consider diffeomorphisms, denoted as $\Phi_\tau^Y$, corresponding to 
the  time--$\tau$ flow of a vector--field $Y$.
Namely, $$z=\Phi^Y_\t(y)=y+Y(y) \t+\cdots$$ is defined as the solution
of
$$\left\{
\begin{array}{lll}
\partial_{\tau} z=Y(z)\\\\
z|_{\tau=0}=y
\end{array}
\right.$$

\noindent
By~\eqref{new vectorfield},  $\Phi_\tau^Y$
transforms a given vector--field $X$ to
$$Z_\tau(y):=J_\tau^Y(y)^{-1}X\big(\Phi^Y_\tau (y)\big)\quad {\rm with} \quad (J_\tau^Y(y))_{hk}:=\partial_{y_k}\big(\Phi^Y_\tau (y)\big)_h \,.$$
We prove the following identity
\begin{align}\label{eq: derivatives}\frac{d^k}{d\tau^k} Z_\tau(y)%\Big|_{\tau_0}
=J_{\tau%_0
}^Y(y)^{-1}{\cal L}^k_Y X\big(\Phi^Y_{\tau%_0
} (y)\big)\end{align}
which immediately will imply the claimed one:
$$ Z(y)=Z_1(y)=\sum_{k=0}^\infty \frac{1}{k!}\frac{d^k}{d\tau^k} Z_\tau(y)\Big|_{0}=\sum_{k=0}^\infty \frac{{\cal L}^k_Y X\big(y\big)}{k!}=e^{{\cal L}_Y}X(y)\,.$$

\noindent
To prove~\eqref{eq: derivatives}, we use the expansion 
\begin{align}\label{exp2}\Phi^Y_\tau (y)=\Phi^Y_{\tau_0} (y)+Y\big(\Phi^Y_{\tau_0} (y)\big)(\tau-\tau_0)+o(\tau-\tau_0)\end{align}
and
\begin{align}\label{inversion}J_\tau^Y(y)=\Big(\id+J_Y\big(\Phi^Y_{\tau_0} (y)\big)(\tau-\tau_0)\Big)J_{\tau_0}^Y(y)+o(\tau-\tau_0)\qquad J_Y(z)_{hk}=\partial_{z_k}Y_h(z)\,.\end{align}
Equation~\eqref{inversion}
gives
\begin{align}\label{exp3}(J_\tau^Y(\eta))^{-1}=(J_{\tau_0}^Y(y))^{-1}\Big(\id-J_Y\big(\Phi^Y_{\tau_0} (y)\big)(\tau-\tau_0)\Big)+o(\tau-\tau_0)\,.\end{align}
While~\eqref{exp2} gives
\begin{eqnarray}\label{expansion}
X\big(\Phi^Y_\tau (y)\big)&=&X\big(\Phi^Y_{\tau_0} (y)+Y\big(\Phi^Y_{\tau_0} (y)\big)(\tau-\tau_0)+o(\tau-\tau_0)\big)\nonumber\\
&=&X\big(\Phi^Y_{\tau_0} (y)\big)+J_X\big(\Phi^Y_{\tau_0} (y)\big) Y\big(\Phi^Y_{\tau_0} (y)\big)(\tau-\tau_0)+o(\tau-\tau_0)
\end{eqnarray}
Collecting~\eqref{exp3} and~\eqref{expansion}, we then find
\begin{eqnarray*}
Z_\tau(y)&=&J_\tau^Y(y)^{-1}X\big(\Phi^Y_\tau (y)\big)\nonumber\\
&=&J_{\tau_0}^Y(y)^{-1}
\nonumber\\
&&
\Big(\id-J_Y\big(\Phi^Y_{\tau_0} (y)\big)(\tau-\tau_0)+o(\tau-\tau_0)\Big)\Big(X\big(\Phi^Y_{\tau_0} (y)\big)+J_X\big(\Phi^Y_{\tau_0} (y)\big) Y\big(\Phi^Y_{\tau_0} (y)\big)(\tau-\tau_0)\Big)\nonumber\\
&&+o(\tau-\tau_0)\nonumber\\
&=&J_{\tau_0}^Y(y)^{-1}X\big(\Phi^Y_{\tau_0} (y)\big)\nonumber\\
&&+J_{\tau_0}^Y(y)^{-1}\Big(J_X\big(\Phi^Y_{\tau_0} (y)\big) Y\big(\Phi^Y_{\tau_0} (y)\big)-J_Y\big(\Phi^Y_{\tau_0} (y)\big) X\big(\Phi^Y_{\tau_0} (y)\big)\Big)(\tau-\tau_0)+o(\tau-\tau_0)
 \end{eqnarray*}
This expansion shows that
$$\frac{d}{d\tau} Z_\tau(y)%\Big|_{\tau_0}
=\frac{d}{d\tau}\Big(
J_{\tau}^Y(y)^{-1}X
\big(\Phi^Y_{\tau} (y)\big)\Big)%\Big|_{\tau_0}
=J_{\tau%_0
}^Y(y)^{-1}{\cal L}_Y X\big(\Phi^Y_{\tau%_0
} (y)\big) $$
By iteration, we have~\eqref{eq: derivatives}.
$ \qquad \square$

\newpage

\addcontentsline{toc}{section}{References}
 \bibliographystyle{plain}
%\bibliography{REFERENCES.bib}
\def\cprime{$'$} \def\cprime{$'$}
\def\cprime{$'$} \def\cprime{$'$}

\end{document}